\documentclass[aip,jmp]{revtex4-1}       
\usepackage{amssymb}
\usepackage{amsmath}
\newcommand{\bea}{\begin{eqnarray}}
\newcommand{\eea}{\end{eqnarray}}
\newcommand{\R}{\mathbf{R}}

\newcommand{\erfc}{\mbox{erfc}}
\newcommand{\tr}{\mbox{tr}}
\newcommand{\sgn}{\mbox{sgn}}
\renewcommand{\bar}{\overline}

\renewcommand{\bar}{\overline}
\newcommand{\dZ}{Z^{\dagger}}
\newtheorem{thm}{Theorem}
\begin{document}
%%%%%%%%%%%%%%%%%%%%%%%%%%%%%%%%%%%%%%%%%%%%%%%%%%%%%%%
\title{The Ginibre evolution in the large-$N$ limit.}
%%%%%%%%%%%%%%%%%%%%%%%%%%%%%%%%%%%%%%%%%%%%%%%%%%%%%%%
\author{Roger Tribe}
\email{tribe@maths.warwick.ac.uk}
\noaffiliation
\affiliation{Department of Mathematics, University of Warwick, Coventry, CV4 AL, UK}
\author{Oleg Zaboronski}
\email{olegz@maths.warwick.ac.uk}
\noaffiliation
\affiliation{Department of Mathematics, University of Warwick, Coventry, CV4 AL, UK}

\date{\today}
\begin{abstract}
We analyse statistics of the real eigenvalues of $gl(N,\mathbf{R})$-valued Brownian motion 
(the 'Ginibre evolution') in the
limit of large $N$. In particular, we
calculate the limiting two-time  correlation function of spin variables associated with real
eigenvalues of the Ginibre evolution. We also show how the formalism of spin variables can
be used to compute the fixed time correlation functions of real eigenvalues discovered originally by Forrester and Nagao \cite{forrester2007eigenvalue} and Borodin and Sinclair\cite{borodin2009ginibre}. 
\end{abstract}
\maketitle
\section{Introduction}
Let $gl(N,\mathbf{R})$ be the linear space consisting of all $N\times N$ real matrices  
equipped with the Euclidean inner product
\bea
\langle A, B\rangle = \tr(AB^T).
\eea
The paper is dedicated to the study of $gl(N,\mathbf{R})$ Brownian motion $(M_t: t\geq 0)$,
the Gaussian process with values in $gl(N,\mathbf{R})$ with zero mean
and covariance given by
\bea
\mathbb{E}\left(\langle A,M_t \rangle \langle B,M_s \rangle \right)
= \frac12 \langle A,B \rangle \, (s\wedge t),~ A,B \in gl(N,\mathbf{R}),
\eea
where $s\wedge t=min(s,t)$ (a standard notation adopted in the probability literature).  
In other words, $(M_t:t \geq 0)$ is the process of $N \times N$ matrices whose entries are 
independent one-dimensional
Brownian motions. The one-dimensional density of $gl(N,\mathbf{R})$ Brownian motion with 
respect to Lebesgue measure is
Gaussian,
\bea
\gamma^{(N)}_t(M) = (\pi t)^{- \frac{N^2}{2}} e^{ - \frac{1}{t} \langle M, M\rangle},
%\label{gpdf}
\eea
and, for $t=1$, this defines the probability measure for the real Ginibre matrix ensemble \cite{mehta2004random}.
We will therefore refer to the process $(M_t, t\geq 0)$ as the real Ginibre evolution, or simply the Ginibre evolution.
The $n$-dimensional Lebesgue density for the Ginibre evolution is
\bea\label{ndimd}
\gamma^{(N)}_{t_1, t_2,\ldots t_n}(M_1, M_2, \ldots, M_n)=\prod_{k=1}^n   
\frac{e^{ - \frac{1}{(t_k-t_{k-1})} \langle{ M_k-M_{k-1}, M_{k}-M_{k-1}\rangle}}}{(\pi (t_k-t_{k-1}))^{\frac{N^2}{2}}},
\eea
where $t_i>t_{j}$ for $i>j$ and $M_0=0,~ t_0=0$.

The principal subject of our study is the stochastic evolution of real eigenvalues of $M_t$ induced by $gl(N,\R)$ 
Brownian motions. 
(For $N>>1$ there are $O(\sqrt{N})$ such eigenvalues \cite{edelman1994many}.)

The closest
counter-part of this process in random matrix theory is the celebrated Dyson Brownian motion defined 
as the eigenvalue process induced by the Brownian motion with values
in Hermitian matrices, see \cite{zeitouni2010}, \cite{mehta2004random} for a review. Of course, the 
nature of the Ginibre evolution is very different: unlike the eigenvalues of Hermitian matrices,  
the real Ginibre eigenvalues can collide and escape into the complex plane and, conversely, a pair of 
complex conjugated eigenvalues can 'land' on the real axis and give birth to the pair
of real eigenvalues. 

Intuitively, it seems conceivable that the large-$N$ evolution of the set of real eigenvalues 
induced by the Ginibre evolution  is described by a Markovian interacting particle 
system on a line, such that particles
are allowed to collide and annihilate and there is a stochastic mechanism  
for the creation of  pairs of particles. 
Indeed, it was shown in \cite{EJP942} that the one dimensional 
distribution of real eigenvalues for the Ginibre evolution
converges as $N \to \infty$ to the one dimensional distribution of 
particles for annihilating Brownian motions on $\R$,
under a suitable initial condition (entrance law). Recall that annihilating Brownian motions is a 
classical interaction particle system. It consists of
particles on a line performing independent diffusions and annihilating on contact. 
Annihilating Brownian motions can
be viewed as the continuous limit of the system of domain walls in the kinetic Glauber spin chain. 
They have been studied by both physics (Glauber, Peilit, Doi, Zeldovich, Ovchinnikov, Derrida, 
Hakim, Pasteur, Lee and Cardy, $\ldots$) and 
mathematics (Bramson and Lebowitz, Griffeath, Kesten, ben-Avraham, Masser $\ldots$) 
communities. See \cite{EJP942} for references and a review of the latest results. 
In order to understand this connection further, one needs to study
the Ginibre evolution beyond the one-dimensional distribution.
As it turns out, such a study can be simplified if carried out in terms of 
'spin' variables associated with real eigenvalues
which we will define as follows. For an $N \times N$ matrix $M$ let $\Lambda^M$ 
be a counting measure on $\R$:
\bea
\Lambda^M(a,b) = \mbox{Number of real eigenvalues of $M$ lying in $(a,b)$.}
\eea
The spin variable associated with $M$ is a function $s(M):\R \rightarrow \{\pm 1\}$:
\bea
s_x(M) = (-1)^{\Lambda^M(-\infty,x)}, \quad x \in \R.
\label{spins}
\eea
Note an analogy between the spin variables (\ref{spins}) and spins in a 
one-dimensional spin chain with real
eigenvalues playing the role of domain walls. Spin variables are crucial in 
linearising the moment equations for
annihilating random walks and/or Brownian motions, see 
e.g. \cite{glauber1963time}, \cite{masser2001method}. 
We believe they will be useful for random matrix models
where eigenvalues are real or the complex eigenvalues appear in pairs. 
Indeed,
the following elementary remark provides a tool for computing 
correlation functions of spin variables: the
spectrum of a real $N\times N$ matrix $M$
consists of real eigenvalues and pairs of conjugated complex eigenvalues. Therefore, 
\bea\label{new7}
s_x(M)=(-1)^{\# \{\mbox{ All eigenvalues of $M$ with real parts in $(-\infty,x)$}\}}
\eea
As a pair of complex conjugated eigenvalues corresponds to a positive 
factor in the characteristic
polynomial, (\ref{new7}) implies that when $M$ has no real eigenvalue at $x$, an event of 
probability $1$ under the Ginibre density,
\bea\label{spin2det}
s_x(M)=\sgn\left(\det\left(M-xI\right)\right).
\eea
We will show that when computed with the help of Householder 
transformations \cite{householder1958unitary}, the correlation functions of spin variables
reduce to the correlation
functions of characteristic polynomials for the Ginibre evolution. 
The latter problem has been thoroughly
investigated, at least in the context of real Ginibre matrix ensembles 
\cite{akemann2008characteristic}, \cite{sommers2009schur}.
All multi-time densities for real eigenvalues can be restored from the multi-time
correlation functions of spin variables. Namely we have the following relation:
\bea
&& \hspace{-.4in} \rho^{(N)}_{t_1, t_2,\ldots, t_K}(x_1, x_2,\ldots, x_K) \nonumber \\
& = & \left. \left(-\frac{1}{2}\right)^K
\left(\prod_{k=1}^K \frac{\partial}{\partial y_{k}}\right)
\mathbb{E}\left(\prod_{m=1}^{K} s_{x_m}\left(M_{t_m}\right)s_{x_m+y_m}\left(M_{t_m}\right) \right)
\right|_{y_m=0+, \, m=1,2\ldots , K}
\label{eq:dens}
\eea
where $\rho^{(N)}_{t_1, t_2,\ldots, t_K}(x_1, x_2,\ldots, x_K)$ is the $K$-dimensional probability density for
finding a real eigenvalue in each of $K$ disjoint intervals at times $t_k$, $k=1,2,\ldots K$
(see \cite{EJP942} for details of (\ref{eq:dens}) in the special case of $t_1=t_2=\ldots =t_K$).
From a purely technical point of view, it is also useful to consider derivatives of $\rho^{(N)}$
leading us to modified densities defined as follows:
\bea
 \tilde{\rho}^{(N)}_{t_1, t_2, \ldots, t_K}(x_1,x_2,\ldots, x_K) \prod_{k=1}^K dx_k
= \mathbb{E} \left( \prod_{k=1}^K s_{x_k}(M_{t_k})  \Lambda^{M_{t_k}}(dx_k) \right).
\label{rhom}
\eea
The above formula is an equality between measures acting on direct products of disjoint intervals (with $dx_k$ on the left hand side
being a standard abbreviation for the Lebesgue measure on $\R$).
Spin correlation functions can be restored from the modified densities
by $n$-dimensional integration, see (\ref{eq:ss1d}). 

In keeping with the standard terminology adopted in probability theory in 
general and random matrix theory in particular, we will 
often refer to densities defined in (\ref{eq:dens}) as correlation functions. 
More precisely, the $n$-dimensional $K$-point correlation function is the 
Lebesgue density for the distribution of $K$ real eigenvalues at $n$ distinct 
time slices. Clearly, $K\geq n$. Therefore, one-dimensional densities refer to equal time
correlation functions , two-dimensional densities - two-time correlation functions and so on.  

The main result of the paper is the exact expression for 
the two-dimensional correlation function of spin variables
 $\mathbb{E}(s_x(M_t)s_{y}(M_{t+\tau}))$. This is done by a lengthy 
calculation, and is
an admittedly modest step towards the complete understanding of the Ginibre evolution.
However, to our knowledge, no multi-time statistics for the Ginibre evolution have ever been calculated.
Moreover this one simple statistic suggests several insights into the general properties of the Ginibre evolution
(see the discussion in the subsequent section). 
In addition, this paper also establishes a novel integral representation for the
fixed time multi-point statistics, which we believe will create a crucial reference point for the future multi-time analysis.

The rest of the paper is organised as follows.  In section \ref{sec:results} 
we state and discuss our main results: the large
$N$ limit of the two-dimensional spin correlation function and a new integral
representation for the one-dimensional multi-point correlation functions of real eigenvalues.
In section \ref{sec:dervtns} we present the main steps of the computation of the 
two-dimensional spin correlation function
for fixed $N$, based on the Householder transform and the technique of integrals 
over anti-commuting variables
(Berezin integrals). In section \ref{largeN} we examine the $N \to \infty$ limit. In section \ref{oned}
we test the techniques developed in section \ref{sec:dervtns} to re-derive the result of
\cite{borodin2009ginibre}, \cite{forrester2007eigenvalue} for the large $N$ limit 
of fixed time multi-point correlation functions
$\rho^{(N)}_{1,\ldots,1}(x_1,\ldots,x_K)$ and $\tilde{\rho}^{(N)}_{1,\ldots,1}(x_1,\ldots,x_K)$, that is the expressions
(\ref{eq:dens}) and (\ref{rhom}) for $t_1=t_2=\ldots=t_K=1$.
%
%%%%%%%%%%%%%%%%%%%%%%%%%%%%%%%%%%%%%%%%%%%%%%%%%%%%%%%%%%%%%%%5
\section{Discussion of results}
\label{sec:results}
Our first result is the de-correlation of spin variables in the real Ginibre evolution for a fixed positive time lag.
\begin{thm}\label{thm1} (Propagation of temporal chaos.) For $x,y \in \mathbf{R}$ and $t,\tau>0$
\bea
\lim_{N\rightarrow \infty} \mathbb{E} \left(s_x(M_{t+\tau})s_y(M_{t})\right) = 0.
\eea
\end{thm}
This shows the large-$N$ and small time lag limits do not commute. 
It also suggests that in the large-$N$
limit, the stochastic evolution of the real eigenvalues associated with $gl(N, \R)$ 
Brownian motion becomes memoryless, that is its marginals at distinct times become independent. 
This result supports our belief that for the real Ginibre evolution
interactions between the $O(\sqrt{N})$ real eigenvalues separated by distances 
of order $1$ are 'screened' by a mean
field created by long-range interactions with the $O(N)$ complex eigenvalues.  
For a long-range first order system of $N$ particles,
the scale of temporal correlations is expected to be of the order of $1/N$,
see \cite{bouchet2004stochastic} for a good heuristic discussion of the Hamiltonian case.

To study short-scale time
correlations we consider the scaling limit $\tau=\frac{T}{N}$ and $N \to \infty$. Therefore
our second result concerns the behaviour of spin-spin correlation function in the real Ginibre evolution
for time lags of order $1/N$.
\begin{thm}\label{thm2} (Decay of temporal correlations.) For $x,y \in \mathbb{R}$ and $t,T>0$.
\bea\label{sscf}
\lim_{N\rightarrow \infty} \mathbb{E} \left(s_x(M_{t+\frac{T}{N}}) \, s_y(M_{t}) \right)
=\erfc\left(\sqrt{\frac{(x-y)^2}{t}+\frac{T}{2t}}\right),
\eea
where $\erfc$ is the complementary error function.
\end{thm}

Notice that for for $T=0$, (\ref{sscf}) reduces to a well-known answer for continuous limit of the
single time spin-spin correlation function in the Glauber model \cite{glauber1963time}. 
This result is not unexpected: as mentioned above, the one-dimensional law of real eigenvalues for
$gl(N)$-Brownian motions coincides with the one-dimensional law for annihilating
Brownian motions on the real line, the latter playing the role of domain walls for
the spin variables \cite{EJP942}. Still, we find this connection very surprising.

As $T\rightarrow \infty$, the spin-spin correlation function decays exponentially, 
which is also very natural if the scenario of 'screening' of real eigenvalues
due to their interactions with complex eigenvalues holds true.

The large-$N$ limit of the two-time spin-spin correlation function characterized by Theorem \ref{thm1}
and Theorem $\ref{thm2}$ is different from the two-time
spin-spin correlation function for the system of one-dimensional 
annihilating Brownian motions obtained in \cite{ECP2133}.
Unfortunately, this means that the conjecture concerning the law of 
real eigenvalues in the real Ginibre evolution made in \cite{EJP942} is false.

On the positive side, Theorem \ref{thm2} yields important 
information about the dynamics of real eigenvalues at small time scales: 
notice that the spin correlation function evaluated at $x=y=0$  
admits the following representation:
\bea\label{spinrep}
\mathbb{E} \left(s_0(M_{t+\tau}) \, s_0(M_{t}) \right)=1-2Pr(\{N_0(t,t+\tau)\mbox{ is odd}\}),
\eea 
where $N_0(t,t+\tau)$ is the number of real eigenvalues crossing $x=0$ 
in the time interval $(t, t+\tau)$. ('Crossing' means that the eigenvalue's 
positions at $t$ and $t+\tau$ straddle the point $x=0$.)
The average density of real eigenvalues  is $1/\sqrt{\pi t}$ and the 
probability of finding two eigenvalues separated by $d<<\sqrt{t}$
is much smaller than the eigenvalue density squared, see e.g. 
Theorem 1 of \cite{EJP942}. Therefore (\ref{spinrep}) suggests
\bea\label{spinrepst}
\mathbb{E} \left(s_0(M_{t+\tau}) \, s_0(M_{t}) \right)=1-2Pr(\{N_0(t,t+\tau)=1\})+o(\tau^{1/2}),
\eea 
Comparing (\ref{spinrepst}) with (\ref{sscf}) we find that
\bea\label{prxing}
Pr(\{N_0(t,t+\tau)=1\})=\sqrt{\frac{N\tau}{2\pi t}}+o(\tau^{1/2})
\eea
(Notice that the rigorous proof of the above relation requires 
expression (\ref{integral}) for the modified density valid at $N<\infty$.)
Recall the small-time behaviour of zero crossing probability in a system of 
rate-$D$ independent diffusive particles scattered over $\R$ with the Poisson intensity $\rho(t)=1/\sqrt{\pi t}$:
\bea\label{prxingm}
Pr_D(\{N_0(t,\tau+\tau)=1\})=\sqrt{\frac{4D\tau}{\pi^2t}}+o(\tau^{1/2})
\eea
Comparing (\ref{prxing}) with (\ref{prxingm}) we conclude that (\ref{sscf}) 
is compatible with the {\it diffusive behaviour of real eigenvalues at small time scales} with an $N$-dependent effective diffusive rate:
\bea
D_{eff}=\frac{N\pi}{8}.
\eea
The linear growth of $D_{eff}$ with $N$ is consistent with the 
de-correlation of multi-time statistics of real eigenvalues in the limit $N\rightarrow \infty$, see Thm. \ref{thm1}.  
Our analysis suggests a possibility of describing the stochastic 
evolution of real eigenvalues at small time scales in terms of a system of 
stochastic differential equations driven by white-in-time noise, even though
the structure of the interaction terms should be very different from 
that of Dyson Brownian motions! 

%However, simulations seem to show that behaviour near collision times between real eigenvalues dominate the
%evolutions, and may be responsible for the de-correlation in time. Note the exponential de-correlation implied by the
%error function is consistent with other long range systems, and with our belief that the large $N$ dynamics may be dominated by
%fluctuations around a mean field limit.

Our final result concerns a novel integral representation for the 
modified $K$-point density of real eigenvalues
at a fixed time.
\begin{thm}\label{thm3} (Ginibre ensemble and anti-self dual Gaussian symplectic ensembles.)
Let $K$ be an even natural number. Fix $x_1<x_2<\ldots <  x_K \in \mathbf{R}$. 
Then
\bea\label{rho_mod}
\lim_{N\rightarrow \infty} \tilde{\rho}^{(N)}_{1,\ldots,1}(x_1, x_2, \ldots, x_K)=C_K\Delta(\mathbf{x})  \int_{U(K)} \mu_H(dU)
e^{-\frac{1}{2}Tr\left(H-H^R\right)^2},
\eea
where $C_K$ is a positive constant, $H=UXU^\dagger$ is a Hermitian 
matrix with eigenvalues $x_1, x_2,\ldots, x_K$, $\mu_H$ is Haar measure
on the unitary group $U(K)$, $H^R=JH^TJ$ is a symplectic involution of 
matrix $H$, $J$ is the canonical symplectic matrix,
and $\Delta(\mathbf{x})=
\prod_{i>j}(x_i-x_j)$ is the Vandermonde determinant.
\end{thm}
In the present paper we only use Theorem \ref{thm3} to show how the formalism 
of spin variables can be employed to re-derive
the Pfaffian expressions for the correlation functions of real eigenvalues in the real 
Ginibre ensemble. These were originally
obtained in \cite{borodin2009ginibre}, \cite{forrester2007eigenvalue}. However we suspect 
that there is a generalization of the Theorem \ref{thm3}
to the case of multi-time correlation functions which should be useful for a 
complete analytic description of the Ginibre evolution.

An additional reason for presenting Theorem \ref{thm3} here is a certain 
mathematical beauty associated with the integral in the right hand side of 
(\ref{rho_mod}). To re-derive Borodin-Sinclaire-Forrester-Nagao Pfaffian 
formula from (\ref{rho_mod}) 
one has to evaluate the integral in the right hand side. The result is (see section V for details):
\bea\label{iz12}
\int_{U(K)} \mu_H(dU)
e^{-\frac{1}{2}Tr\left(H-H^R\right)^2}=C'_K  \frac{Pf\left[(x_i-x_j)e^{-2(x_i-x_j)^2}\right]_{1\leq i,j\leq K}}{\Delta(\mathbf{x})},
\eea
where $C'_K$ is a positive constant. Intriguingly, the right hand side of the 
above formula is the first term of the stationary phase
expansion of the integral on the left hand side.\footnote{As conjectured by 
Yan Fyodorov during an after-seminar discussion.} We verified this in complete analogy with the proof
of the exactness of the stationary phase expansion for the Itzykson-Zuber integral 
carried out in \cite{stone}. In other words, we checked that the integral in 
the left hand side of (\ref{iz12}) localizes exactly to the set of stationary points of $F(H)=Tr\left(H-H^R\right)^2$.
At the moment, the precise reason for this localization is unclear to us. 
In particular, the Duistermaat-Heckmann theorem\cite{getzler} which is 
responsible for the exact localization of the Itzykson-Zuber-Harish-Chandra
integral is not directly applicable to our case. Due to symplectic invariance 
of $F(H)$, the integral in  (\ref{iz12}) is taken over the symmetric 
space $U(K)/USp(K)$, where $USp(K)$ is the symplectic subgroup of $U(K)$. 
But $dimU(K)/USp(K)=K(K-1)/2$, which is even only if $K$ is divisible 
by $4$. So in general,   $U(K)/USp(K)$ is not even symplectic and the Duistermaat Heckmann theorem does not apply. 

 A variant of the integral (\ref{iz12}) plays an important role in the 
analysis of multi-dimensional multi-point correlation functions for the real 
Ginibre evolution\cite{tribe_fc}. We feel that a proper understanding
of the localization of (\ref{iz12}) will lead to a significant further progress in 
uncovering the structure of this stochastic process.

Theorems \ref{thm1} and \ref{thm2} are proved in sections \ref{sec:dervtns}, \ref{largeN}. Theorem
\ref{thm3} is proved in section \ref{oned}.

%
%%%%%%%%%%%%%%%%%%%%%%%%%%%%%%%%%%%%%%%%%%%%%%%%%%%%%%%%%%%%%
%%%%%%%%%%%%%%%%%%%%%%%%%%%%%%%%%%%%%%%%%%%%%%%%%%%%%%%%%%%%%%%%%
\section{Calculation of the two-dimensional spin-spin correlation function.}
\label{sec:dervtns}
%%%%%%%%%%%%%%%%%%%%%%%%%%%%%%%%%%%%%%%%%%%%%%%%%%%%%%

The calculation detailed below draws on the ideas of Edelman et. al. \cite{edelman1994many} and Sommers et. al.
\cite{sommers2008general} to compute the
two-dimensional spin-spin correlation function while avoiding the mathematically
difficult transition to the eigenvalue representation. 
%Edelman used a Householder transform
%to compute the density of the eigenvalues in the real Ginibre ensemble while Sommers
%employed Berezin integrals to simplify his derivation of the kernel of the Pfaffian point
%process associated with the real Ginibre ensemble.
%while avoiding the mathematically difficult transition to the eigenvalue representation. 
Edelman used a
Householder transform to compute the density of the eigenvalues in the real Ginibre ensemble
while Sommers employed Berezin integrals to simplify his derivation of the kernel of the Pfaffian
point process associated with the real Ginibre ensemble.
The key fortuitous link which makes our calculation possible is that $s_x(A) = \sgn( \det(A-xI))$
(see eq. (\ref{spins})) which combines nicely with the Jacobian
$|\det(A-xI)|$ of the
Edelman transformation (see below) to produce the characteristic polynomial $\det(A-xI)$. The latter is easier
to average over the Ginibre evolution.

Instead of computing the spin-spin correlation function directly, we
compute the time-space point modified density function $\tilde{\rho}^{(N)}$ first:
\bea
 \tilde{\rho}^{(N)}_{t,t+\tau}(y,x)dxdy
= \mathbb{E} \left[ s_y(M_t)  \Lambda^{M_t}(dy)  s_x(M_{t + \tau}) \Lambda^{M_{t + \tau}}(dx) \right].
\label{rhot}
\eea
The above notation stands stands for equality of measures on $\R^2$. 
To recover the density $\tilde{\rho}^{(N)}$ it is crucial to take infinitesimal intervals $dx$ and $dy$
to lie just to the right of the points $x$ and $y$ respectively. For example, the measure
$\mathbb{E} [s_a(M_t) \Lambda^{M_t}(db)]$ has a density which is discontinuous at $a$,
indeed the sign of the density flips.
Given $\tilde{\rho}^{(N)}$, the spin-spin correlation function can be restored by integration:
\bea\label{inv2}
\mathbb{E}\left(s_x(M_{t+\tau})s_y(M_t)\right)= 4
\int_{-\infty}^x \!\! dx' \int_{-\infty}^y \!\! dy' \,  \tilde{\rho}^{(N)}_{t,t+\tau}(y',x').
\eea
It follows from the definition (\ref{rhot}) of the modified density that
\begin{eqnarray}
& & \hspace{-.3in} \tilde{\rho}^{(N)}_{t,t+\tau} (y,x)dxdy \nonumber \\
& = & \int_{R^{N^2}} dM \int_{R^{N^2}} dH  s_y(H) \Lambda^H(dy)  s_x(M)  \Lambda^M(dx)
\, \gamma_t(H) \, \gamma_{\tau}(M-H). \label{zero}
\end{eqnarray}
The plan of attack is, via changes of variable that exploit symmetry of the integrand,
to reduce the dimension of this integral. Indeed by the end of this section the
integral over $\R^{2N^2}$ is reduced to a one dimensional integral over $[0,1]$, 
see (\ref{integral}).
\noindent
\textbf{The representation of $dM$ integral in Edelman coordinates.}
The transformation used by Edelman et. al. in \cite{edelman1994many} is for an $N \times N$ matrix $M$
with a real eigenvalue $x$ and corresponding eigenvector
$v \in S^+_{N-1}$, the upper half of the $N-1$ dimensional unit sphere in $\R^N$:
\bea\label{edt}
M = P_v M^e P_v
\eea
where $P_v$ is the Householder transformation \cite{householder1958unitary} 
that reflects in the hyperplane at right
angles to the vector $v-e_N$ (where $e_N$ is the unit vector $(0,\ldots,0,1)$), 
and $M^e$ is in Edelman block form
\bea\label{edt1}
M^e = \left( \begin{array}{cc}
M^e_0 & 0 \\
w^T & x
\end{array} \right)
\eea
with $M_0^e$ an $(N-1) \times (N-1)$ real matrix, $w \in R^{N-1}$ and $x \in \R$. 
Later we need the explicit form
\begin{equation} \label{pvform}
P_v = I - \frac{(v-e_N)(v-e_N)^T}{1-v_N}
\end{equation}
where $v_N<1$ is the $N$th co-ordinate of $v$.  Note that $P_v$ is orthogonal 
and symmetric, $P_v^{-1}=P_v^T = P_v$.
Let
\bea
F(M) dy = \int_{R^{N^2}} dH s_y(H) \Lambda^H(dy) \gamma_t(H) \gamma_{\tau}(M-H).
\eea
Fix $t,\tau>0$ and $x,y$ throughout and write $\tilde{\rho}^{(N)}$ as 
shorthand for $\tilde{\rho}^{(N)}_{t,t+\tau} (y,x)$.
Using the known expression for the Jacobian of the Edelman transformation 
\cite{edelman1994many}, we rewrite
the expression (\ref{zero}) for $\tilde{\rho}^{(N)}$ as
\begin{eqnarray}
\tilde{\rho}^{(N)} & = & \int_{S_{N-1}^+} \! \! dv \int_{R^{(N-1)^2}} dM^e_0 \int_{R^{N-1}} dw
|\det(M_0^e-xI)| s_x(M_0^e) F(P_v M^e P_v)\nonumber \\
& = & \int_{S_{N-1}^+} \! \! dv \int_{R^{(N-1)^2}} dM^e_0 \int_{R^{N-1}} dw
 \det(M_0^e-xI) F(P_v M^e P_v)\label{abs_val}
\end{eqnarray}
where $dv$ stands for surface measure on the sphere. We have 
also used $s_x(P_v M^eP_v) = s_x(M_e)$ since
$s_x(A)$ depends only on the eigenvalues of $A$, and
$s_x(M^e) = s_x(M^e_0)$.  To obtain the second equality in (\ref{abs_val}) we used expression
(\ref{spin2det}) relating the value of the spin at $x$ to the sign of the characteristic polynomial
$\det(M-xI)$. Note the disappearance of the absolute value sign from the integrand
in the last expression: the problem of computing the spin-spin correlation
function has been reduced to the problem of averaging the characteristic
polynomial over $M_0^e$.

\noindent
\textbf{Evaluation of the $dv$ integral.} %OK
%%%%%%%%%%%%%%%%%%%%%%%%%%%%%%%%%%%%%%%%%%%%%%%%%%%%%%%%
We claim that $F(O^TMO) = F(M)$ for an orthogonal matrix $O$. To see the latter note that the
Gaussian density has this invariance and therefore
\begin{eqnarray}
F(O^TMO) dy & = & \int_{R^{N^2}} dH s_y(H) \Lambda^H(dy) \gamma_t(H) \gamma_{\tau}(O^TMO-H) \nonumber\\
& = & \int_{R^{N^2}} dH s_y(H) \Lambda^H(dy) \gamma_t(H) \gamma_{\tau}(M- OHO^T) \nonumber\\
& = & \int_{R^{N^2}} dH s_y(O^THO) \Lambda^{O^THO}(dy) \gamma_t(H) \gamma_{\tau}(M-H) \nonumber\\
& = & F(M) dy.
\end{eqnarray}
Here we used that $s_y(H)$ and $\Lambda^{H}(dy)$ depend only on the eigenvalues of $H$ and are therefore
invariant with respect to the transformation $H \to O^THO$.
Then
\begin{eqnarray}
\tilde{\rho}^{(N)} & = & \int_{S_{N-1}^+} dv \int_{R^{(N-1)^2}} dM^e_0 \int_{R^{N-1}} dw
 \det(M_0^e-xI)  F(M^e) \nonumber \\
& = & \frac12 |S_{N-1}|  \int_{R^{(N-1)^2}} dM^e_0 \int_{R^{N-1}} dw
 \det(M_0^e-xI)  F(M^e) \label{basepoint1}
\end{eqnarray}
since the integrand is independent of $v$.  Here $|S_{N-1}|$ stands for the surface area of the unit
$(N-1)$ dimensional sphere living in $\R^N$.

\noindent
\textbf{Evaluation of the $dw$ integral.} %OK
%%%%%%%%%%%%%%%%%%%%%%%%%%%%%%%%%%%%%%%%%%%%%%%%%%%%%%%%%%%
To integrate over $w$ we must express $F(M^e)$ in terms of $M^e_0$ and $w$. Let us represent $H$ in
block form:
for $z \in \R$ and $\alpha,\beta \in \R^{N-1}$, write
\begin{equation} \label{H1stform}
H = \left( \begin{array}{cc}
H_0 & \beta \\
\alpha^T & z
\end{array} \right).
\end{equation}
Expanding the Gaussian densities we find that $F(M^e) dy$ is given by
\begin{equation}\label{eq:F}
 (\pi^2 t \tau)^{- \frac{N^2}{2}} \int_{R^{N^2}} dH s_y(H) \Lambda^H(dy) \;
e^{- \left(\frac{1}{t}+ \frac{1}{\tau}\right) \tr(H H^T)} \; e^{ - \frac{1}{\tau} \tr(M^e M^{eT})} \;
e^{ \frac{2}{\tau} \tr(M^e H^T)}.
\end{equation}
The traces expand in block form to
\begin{eqnarray}
\tr(M^e M^{eT}) & = & \tr(M^e_0 M_0^{eT}) + |w|^2 + x^2, \nonumber \\
\tr(M^e H^T) & = & \tr(M_0^e H_0^T) + w^T \alpha + xz. \label{traces}
\end{eqnarray}
Substituting (\ref{eq:F}) and (\ref{traces}) into (\ref{basepoint1})
we find the following representation for $\tilde{\rho}^{(N)}$:
\begin{eqnarray}
\tilde{\rho}^{(N)} dy & = & \frac12 |S_{N-1}|   (\pi^2 t \tau)^{- \frac{N^2}{2}}
e^{-\frac{x^2}{\tau}}  \int_{R^{(N-1)^2}} dM^e_0 \int_{R^{N-1}} dw
\int_{R^{N^2}} dH   \nonumber \\
&& \hspace{.2in} s_y(H) \Lambda^H(dy) \det(M_0^e-xI)
\; e^{- \left(\frac{1}{t}+ \frac{1}{\tau}\right) \tr(H H^T)} \;
 \nonumber \\
 && \hspace{.4in}
 e^{-\frac{1}{\tau} (|w|^2 - 2 w^T \alpha - 2xz)} \; e^{ - \frac{1}{\tau} \tr(M^e_0 M^{eT}_0)}
 \; e^{ \frac{2}{\tau} \tr(M^e_0 H_0^T)}. \label{basepoint2}
\end{eqnarray}
The $dw$ integral is then computed using a standard formula for Gaussian integrals:
\begin{equation} \label{eq:dw}
\int_{R^{N-1}} dw \,  e^{-\frac{1}{\tau} (|w|^2 - 2 w^T \alpha)} = (\pi \tau)^{\frac{N-1}{2}} e^{\frac{1}{\tau} |\alpha|^2}.
\end{equation}

\noindent
\textbf{Evaluation of the $dM^e_0$ integral.} %OK
%%%%%%%%%%%%%%%%%%%%%%%%%%%%%%%%%%%%%%%%%%%%%%%%%%%
To compute the $dM^e_0$ integral we first re-express the determinant as a Berezin integral over
anti-commuting (Grassmann) variables:
\bea\label{gbi}
\det(M_0^e-xI) = \int_{R^{0|2(N-1)}} d \phi d \bar{\phi} \; e^{\bar{\phi}^T (M^e_0 - xI) \phi}.
\eea
This leaves the following  $dM^e_0$ integral in (\ref{basepoint2})
\begin{equation} \label{I1}
I_1 := \int_{R^{(N-1)^2}} dM^e_0 \; e^{-\frac{1}{\tau}\tr(M_0^e M_0^{eT})} \; e^{\frac{2}{\tau}  \tr(M^e_0 H_0^T) }
\; e^{ \bar{\phi}^T M^e_0 \phi}.
\end{equation}
This is a Gaussian integral with identity covariance matrix
but with the linear term of the exponent depending on anti-commuting variables.
Nevertheless the rules for this integral are as if it were a standard Gaussian
(see \cite{itzykson1991statistical} for details), and the value can be found by
locating the critical value of the quadratic form in the variables $M^e_0(i,j)$.
The integrand in  (\ref{I1}) is the exponential of
\begin{equation} \label{temp143}
- \frac{1}{\tau} \tr(M^e_0 M^{eT}_0) + \frac{2}{\tau} \tr(M^e_0 \Gamma^T)
\end{equation}
where $\Gamma = H_0 + \frac{\tau}{2} \bar{\phi} \phi^T$.
The critical point of this quadratic form is $M_0^e = \Gamma$.
Substituting this back into (\ref{temp143}) one obtains, noting that $(\bar{\phi} \phi^T)^T = - \phi \bar{\phi}^T$,
\begin{eqnarray}
+ \frac{1}{\tau} \tr(\Gamma \Gamma^T)
&=&  + \frac{1}{\tau} \tr((H_0 + \frac{\tau}{2} \bar{\phi} \phi^T)(H_0^T - \frac{\tau}{2} \phi \bar{\phi}^T)) \nonumber\\
& = & + \frac{1}{\tau} \tr(H_0 H_0^T) + \tr(H_0 \bar{\phi} \phi^T)
\end{eqnarray}
This gives the value
\begin{equation} \label{appendix1}
I_1 = (\pi \tau)^{\frac{(N-1)^2}{2}} e^{\frac{1}{\tau} \tr(H_0 H_0^T) + \bar{\phi}^T H_0 \phi}.
\end{equation}
Substituting in the $dw$ and $dM^e_0$ integrals (\ref{eq:dw}) and (\ref{appendix1})
into (\ref{basepoint2}) we reach
\begin{eqnarray}
\tilde{\rho}^{(N)} dy & = &
\frac12 |S_{N-1}|  (\pi t)^{-\frac{N^2}{2}}  (\pi \tau)^{-\frac{N}{2}} e^{-\frac{x^2}{\tau}} \int_{R^{N^2}}
dH \int_{R^{0|2(N-1)}} d \phi d\bar{\phi} \nonumber \\
& &  \hspace{.05in}  s_y(H) \Lambda^H(dy) \; e^{- \left(\frac{1}{t}+ \frac{1}{\tau}\right) \tr(H H^T)} \;
 e^{\frac{1}{\tau} (2xz + |\alpha|^2)} \;  e^{\frac{1}{\tau} \tr(H_0 H_0^T) + \bar{\phi}^T (H_0 -xI) \phi}.
\end{eqnarray}
The integral at hand can be simplified by an orthogonal transformation on the $H$ variables.
To simplify the implementation of the transformation, we will re-write the above expression
in terms of matrix $H$ rather than its sub-matrix $H_0$.
To this end, we use an extended set of $N$-dimensional Grassmann variables  $\bar{\psi},\psi$ that agree with
$\bar{\phi},\phi$ in the first $(N-1)$ co-ordinates,
\bea
\psi=\{\phi, \psi_N\},~\bar{\psi}=\{\bar{\phi}, \bar{\psi}_N\}.
\eea
Then
\bea
\bar{\psi}^T (H-xI) \psi = \bar{\phi}^T (H_0-xI) \phi + (\bar{\phi}^T \beta) \psi_N + \bar{\psi}_N (\alpha^T \phi)
+ (z-x) \bar{\psi}_N \psi_N,
\eea
and we can re-write
\bea
\int_{R^{0|2(N-1)}} d \phi d\bar{\phi} \; e^{\bar{\phi}^T(H_0-xI)\phi}
= \int_{R^{0|2N}} d \psi d\bar{\psi} \; \bar{\psi}_N \psi_N  \, e^{\bar{\psi}^T(H-xI)\psi}
\eea
since the integral over the pair $d \psi_N d \bar{\psi}_N$ can be done first on the right hand side and
the term $\bar{\psi}_N \psi_N$ forces it to take value one. Using also
$ \tr(HH^T) =  \tr(H_0 H_0^T) + |\beta|^2 + |\alpha|^2 + z^2 $
this leaves
\begin{eqnarray}
\tilde{\rho}^{(N)} dy & = & \frac12 |S_{N-1}|  (\pi t)^{-\frac{N^2}{2}} (\pi \tau)^{-\frac{N}{2}}  e^{-\frac{x^2}{\tau}}
\int_{R^{N^2}} dH \, s_y(H) \Lambda^H(dy)
\nonumber \\
& &  \hspace{.2in}  \int_{R^{0|2N}} d \psi d\bar{\psi} \; \bar{\psi}_N \psi_N \;
e^{\bar{\psi}^T(H-xI)\psi} \; e^{- \frac{1}{t}\tr(H H^T)}
\; e^{- \frac{1}{\tau} (|\beta|^2+ z^2 - 2xz)}.
\label{basepoint3}
\end{eqnarray}

\noindent
\textbf{Representation of the $dH$ integral in Edelman variables.}
%%%%%%%%%%%%%%%%%%%%%%%%%%%%%%%%%%%%%%%%%%%%%%%%%%%%%%
A second Edelman change of variable is $H=P_v H^e P_v$ where
\bea
H^e = \left( \begin{array}{cc}
H^e_0 & 0 \\
w^T & y
\end{array} \right),
\eea
see (\ref{edt}), (\ref{edt1}), (\ref{pvform}) for the full definition of Edelman
transform. 
We need to reconcile this with our earlier representation (\ref{H1stform}) of $H$ 
in block form. Let  $v=(\hat{v},v_N)$, where
$\hat{v}$ is an $(N-1)$-dimensional vector and $v_N$ - a scalar such that
$\hat{v}\cdot \hat{v}+v_N^2=1$.
Multiplying the three matrices $P_v ,~H^e$ and $P_v$ explicitly using (\ref{pvform}) and 
comparing the result with (\ref{H1stform}) we find:
\begin{eqnarray}
z & = & (\hat{v}^TH_0^e \hat{v}) + v_N (w^T \hat{v}) + y v_N^2, \nonumber \\
 \beta & = & H_0^e \hat{v} - (\hat{v}^T H_0^e \hat{v}) (1-v_N)^{-1} \hat{v}
+ (w^T \hat{v}) \hat{v} + y v_N \hat{v}, \nonumber \\
|\beta|^2 + z^2 &=&  (v_Ny + w^T \hat{v})^2 + \hat{v}^T H_0^{eT}H_0^e \hat{v}.
\label{appendix2}
\end{eqnarray}
Next
\bea
\bar{\psi}^T(H-xI) \psi  = (P_v \bar{\psi})^T(H^e-xI) P_v \psi,
\eea
which suggests a change of Grassmann variables $\phi = P_v \psi, \, \bar{\phi} = P_v \bar{\psi}$. Under this
change we find with the help of (\ref{pvform})  that $\psi_N = v^T \phi$ and $\bar{\psi}_N = v^T \bar{\phi}$.
Therefore, Berezin integral in the right hand side of (\ref{basepoint3}) transforms as follows:
\begin{equation}\label{eq:IF}
I_2 : = \int_{R^{0|2N}} d \psi d\bar{\psi} \; \bar{\psi}_N \psi_N \, e^{\bar{\psi}^T(H-xI) \psi}
= \int_{R^{0|2N}} d \phi d\bar{\phi} \, (v^T\bar{\phi}) (v^T\phi) \, e^{\bar{\phi}^T(H^e-xI) \phi}.
\end{equation}
Using this and the Edelman substitution in (\ref{basepoint3}) we reach
\begin{eqnarray}
\tilde{\rho}^{(N)} & = & \frac12 |S_{N-1}|  (\pi t)^{-\frac{N^2}{2}} (\pi \tau)^{-\frac{N}{2}}  e^{-\frac{x^2}{\tau}}
\int_{S_{N-1}^+} \! \! dv \int_{R^{(N-1)^2}} dH^e_0 \int_{R^{N-1}} dw \nonumber\\
&&   \hspace{.2in}   \int_{R^{0|2N}}  d \phi d\bar{\phi} \, (v^T\bar{\phi}) (v^T\phi) \, e^{\bar{\phi}^T(H^e-xI)\phi} \det(H^e_0 - yI) \\
&&   \hspace{.4in} e^{- \frac{1}{t}\tr(H^e H^{eT})} \;
 e^{- \frac{1}{\tau} ((v_Ny + w^T \hat{v})^2 + \hat{v}^T H_0^{eT}H_0^e \hat{v})} \;
e^{ \frac{2x}{\tau}(\hat{v}^TH_0^e \hat{v} + v_N (w^T \hat{v}) + y v_N^2)}.\nonumber
\end{eqnarray}
To prepare for the integration over $w$ and $H^e_0$ we found it convenient to integrate
out some Grassmann variables. First
we let $\hat{\phi},\bar{\hat{\phi}}$ be the
restriction of $\phi,\bar{\phi}$ to the first $(N-1)$ co-ordinates so that
\bea
\bar{\phi}^T (H^e-xI) \phi =  \bar{\hat{\phi}}^T (H_0^e-xI) \hat{\phi}
+ \bar{\phi}_N (w^T \hat{\phi}) + (y-x) \bar{\phi}_N \phi_N.
\eea
We will integrate out the final pair of variables $\phi_N, \bar{\phi}_N$. Extracting the terms
depending on $\phi_N, \bar{\phi}_N$ gives
\begin{eqnarray}
&& \hspace{-.1in} \int_{R^{0|2}} d \phi_N d \bar{\phi}_N
(\hat{v}^T \bar{\hat{\phi}}+ v_N \bar{\phi}_N) (\hat{v}^T \hat{\phi} + v_N \phi_N)
e^{\bar{\phi}_N (w^T \hat{\phi}) + (y-x) \bar{\phi}_N \phi_N} \nonumber\\
& = & \int_{R^{0|2}} d \phi_N d \bar{\phi}_N
(\hat{v}^T \bar{\hat{\phi}}+ v_N \bar{\phi}_N) (\hat{v}^T \hat{\phi} + v_N \phi_N)
(1+ \bar{\phi}_N (w^T \hat{\phi}))(1+ (y-x) \bar{\phi}_N \phi_N) \nonumber\\
& = & v_N^2 + (y-x)(\hat{v}^T \bar{\hat{\phi}}) (\hat{v}^T \hat{\phi}) + v_N (w^T \hat{\phi})
(\hat{v}^T \bar{\hat{\phi}}).
\end{eqnarray}
In the above calculation we used that $ (w^T \hat{\phi})$ and $(\hat{v}^T \bar{\hat{\phi}})$ are odd hence 
nilpotent elements of the Grassmann algebra.
This leaves us just with the integrand depending on shortened anti-commuting 
variables $\hat{\phi}, \bar{\hat{\phi}}$ and
so we can drop the hats in the notation for integration variables replacing 
them with $\phi,\bar{\phi} \in \R^{0|(N-1)}$.
Thus we have shown that the Berezin
integral (\ref{eq:IF}) has become
\begin{eqnarray}
I_2 & = &  \int_{R^{0|2(N-1)}} d \phi d\bar{\phi}
\left(v_N^2 + (y-x)(\hat{v}^T \bar{\phi}) (\hat{v}^T \phi) + v_N (w^T \phi)
(\hat{v}^T \bar{\phi}) \right)   e^{\bar{\phi}^T(H_0^e-xI) \phi} \nonumber\\
& = &  \int_{R^{0|2(N-1)}} d \phi d\bar{\phi}
\left(v_N^2 + (y-x)(\hat{v}^T \bar{\phi}) (\hat{v}^T \phi) \right)
e^{\frac{1}{v_N} (w^T \phi) (\hat{v}^T \bar{\phi})} \;
 e^{\bar{\phi}^T(H_0^e-xI) \phi}.
\end{eqnarray}
Using this, representing $\det(H_0^e - yI)$ as a Gaussian integral over a second set of Grassmann variables
$\bar{\psi}, \psi$ and expanding
\bea
 \tr(H^eH^{eT}) =  \tr(H^E_0 H_0^{eT}) + |w|^2 + y^2
\eea
we reach the following expression for $\tilde{\rho}^{(N)}$:
\begin{eqnarray}
\tilde{\rho}^{(N)} & = & \frac12 |S_{N-1}|  (\pi t)^{-\frac{N^2}{2}} (\pi \tau)^{-\frac{N}{2}}  e^{-\frac{x^2}{\tau}}
\int_{S_{N-1}^+} \! \! dv \int_{R^{(N-1)^2}} dH^e_0 \int_{R^{N-1}} dw \nonumber \\
&&   \hspace{.1in} \int_{R^{0|4(N-1)}}  d \phi d\bar{\phi} d \psi d\bar{\psi}
\, e^{\bar{\phi}^T(H^e_0-xI)\phi + \bar{\psi}^T(H^e_0-yI)\psi} \nonumber \\
&&   \hspace{.2in} \left(v_N^2 + (y-x)(\hat{v}^T \bar{\phi}) (\hat{v}^T \phi) \right)
e^{\frac{1}{v_N} (w^T \phi) (\hat{v}^T \bar{\phi})} \;
 e^{- \frac{1}{t}(\tr(H^e_0 H_0^{eT})+ |w|^2+y^2)} \nonumber \\
&&    \hspace{.3in} e^{- \frac{1}{\tau} ((v_Ny + w^T \hat{v})^2 + \hat{v}^T H_0^{eT}H_0^e \hat{v})} \;
e^{ \frac{2x}{\tau}(\hat{v}^TH_0^e \hat{v} + v_N (w^T \hat{v}) + y v_N^2)}.
\label{basepoint4}
\end{eqnarray}

\noindent
\textbf{Evaluation of the second $dw$ integral.}
%%%%%%%%%%%%%%%%%%%%%%%%%%%%%%%%%%%%%%%%%%%%%%
Collecting all the terms in the integrand of (\ref{basepoint4}) 
depending $w$ one reaches the second $dw$ integral,
which happens to be Gaussian:
\begin{equation} \label{I2}
I_3 = \int_{R^{N-1}} dw e^{- \frac{1}{t}|w|^2} \;
 e^{- \frac{1}{\tau} (2(y-x)v_N (w^T\hat{v}) + (w^T \hat{v})^2)} \;
e^{\frac{1}{v_N} (w^T \phi) (\hat{v}^T \bar{\phi})} .
\end{equation}
The exponential has a quadratic term $t^{-1} w^T M w$ with covariance matrix
\begin{equation}\label{eq:M}
M = I + \frac{t}{\tau} \hat{v} \hat{v}^T.
\end{equation}
Note that
\begin{equation} \label{minverse}
M^{-1} = I - \frac{t}{\tau + t|\hat{v}|^2} \hat{v} \hat{v}^T = I - \frac{t}{\tau D} \hat{v} \hat{v}^T
\quad \mbox{where $D : = \det(M) = 1 + \frac{t}{\tau} |\hat{v}|^2$.}
\end{equation}
The integrand in  (\ref{I2}) is the exponential of
\begin{equation} \label{temp144}
- \frac{1}{t} w^T M w + \frac{2}{t} w^T \gamma
\end{equation}
where $\gamma = \frac{t}{\tau}(x-y)v_N \hat{v} - \frac{t}{2v_N} ( \hat{v}^T \bar{\phi}) \phi$.
The critical value of $w$ in then $M^{-1} \gamma$ and substituting into (\ref{temp144})
one obtains $t^{-1} \gamma^T M^{-1} \gamma$ which is
\begin{eqnarray}
&& \hspace{-.2in}
\frac{1}{t} ( \frac{t}{\tau}(x-y)v_N \hat{v}^T - \frac{t}{2v_N} (\hat{v}^T \bar{\phi}) \phi^T) M^{-1}
(\frac{t}{\tau}(x-y)v_N \hat{v} - \frac{t}{2 v_N} (\hat{v}^T \bar{\phi}) \phi) \nonumber\\
& = & \frac{t}{\tau^2} (x-y)^2 v_N^2 \hat{v}^T M^{-1} \hat{v}
- \frac{t}{\tau} (x-y) (\hat{v}^T \bar{\phi}) \hat{v}^T M^{-1} \phi.
\end{eqnarray}
Using the form (\ref{minverse}) for $M^{-1}$ one finds that $\hat{v}^T M^{-1} = D^{-1} \hat{v}^T$
and this leads to the answer
\begin{equation} \label{appendix3}
I_3 = D^{-\frac12} (\pi t)^{\frac{N-1}{2}} e^{-\frac{t}{\tau D} (x-y) (\hat{v}^T \bar{\phi})(\hat{v}^T \phi)}
\; e^{\frac{t}{\tau^2 D} (x-y)^2 v_N^2 |\hat{v}|^2}.
\end{equation}

\noindent
\textbf{Evaluation of the $dH^e_0$ integral.}
%%%%%%%%%%%%%%%%%%%%%%%%%%%%%%%%%%%%%%%%%%%%%%
Using expression (\ref{eq:M}) for matrix $M$ we can verify the formula
\bea
\frac{1}{t} \tr(H^e_0 H^{eT}_0) + \frac{1}{\tau}  \hat{v}^T H_0^{eT}H_0^e \hat{v} = \frac{1}{t} \tr(H^e_0 M H^{eT}_0).
\eea
Combining all the $H^e_0$-dependent terms in the integrand 
of (\ref{basepoint4}) and using the above identity
we reach the integral
\begin{equation}
I_4 = \int_{R^{(N-1)^2}} dH^e_0 \, e^{\bar{\phi}^T H_0^e \phi + \bar{\psi}^T H^e_0 \psi}
\; e^{- \frac{1}{t} \tr(H^e_0 M H^{eT}_0)}  \; e^{\frac{2x}{\tau} \hat{v}^T H_0^e \hat{v} }.
 \label{I3}
\end{equation}
The integral $I_4$ is Gaussian and can be computed in the usual way by evaluating
the integrand at the critical point of the exponent.
The quadratic part of the exponential $\tr(H^e_0 M H^{eT}_0)$ 
can be considered as a quadratic form in the
variables $(H_0^e(i,j):i,j=1,\ldots,N-1)$. The corresponding 
covariance matrix has determinant
$D^{N-1}$. The easiest way to see this is to take $v = (1,0,\ldots,0)$ and 
examine the effect of the
transformation $H_0^e(i,j) \to (H^e_0M)(i,j)$ for each $i,j$.
The integrand in (\ref{I3}) is the exponential of
\begin{equation} \label{temp145}
- \frac{1}{t} \tr(H^e_0 M H^{eT}_0) + \frac{2}{t} \tr (H^e_0 \Gamma^T)
\end{equation}
 where
\bea
\Gamma = \frac{t}{2} \bar{\phi} \phi^T + \frac{t}{2} \bar{\psi} \psi^T
+ \frac{tx}{\tau} \hat{v} \hat{v}^T.
\eea
The critical point of (\ref{temp145})  is $H^e_0 = \Gamma M^{-1}$.
Substituting into (\ref{temp145})
one obtains $t^{-1} \tr( \Gamma M^{-1} \Gamma^T)$ which is
\bea
\frac{tx^2}{\tau^2} \tr( \hat{v} \hat{v}^T M^{-1} \hat{v} \hat{v}^T)
- \frac{t}{2} \tr ( \bar{\phi} \phi^T M^{-1} \psi \bar{\psi}^T)
+ \frac{tx}{\tau} \tr( \bar{\phi} \phi^T M^{-1} \hat{v} \hat{v}^T)
+  \frac{tx}{\tau} \tr( \bar{\psi} \psi^T M^{-1} \hat{v} \hat{v}^T)\nonumber\\
\eea
and evaluates,  using the form (\ref{minverse}) for $M^{-1}$, to
\bea
\frac{tx^2}{\tau^2 D} |\hat{v}|^4
- \frac{t}{2} (\bar{\phi}^T \bar{\psi}) (\phi^T \psi)
+ \frac{t^2}{2 \tau D} (\bar{\phi}^T \bar{\psi}) (\phi^T \hat{v}) (\psi^T \hat{v})
\nonumber\\ 
+ \frac{tx}{\tau D} (\bar{\phi}^T \hat{v})(\phi^T \hat{v})
 + \frac{tx}{\tau D} (\bar{\psi}^T \hat{v})(\psi^T \hat{v}).
\eea
This gives
\begin{eqnarray}
I_4 &=&  D^{-(N-1)/2} (\pi t)^{(N-1)^2/2}
e^{ - \frac{t}{2} (\phi^T \psi) (\bar{\phi}^T \bar{\psi}) } \;
e^{ \frac{t^2}{2\tau D} (\hat{v}^T \phi) (\hat{v}^T \psi) (\bar{\phi}^T \bar{\psi})}
\nonumber \\
&& \hspace{.2in} e^{  \frac{tx}{\tau D}((\hat{v}^T \bar{\phi})(\hat{v}^T \phi) +
(\hat{v}^T \bar{\psi})(\hat{v}^T \psi))}
\; e^{\frac{tx^2}{\tau^2 D} |\hat{v}|^4}.
\label{appendix4}
\end{eqnarray}
Substituting (\ref{appendix3}) and (\ref{appendix4}) into (\ref{basepoint4}) we reach
\begin{eqnarray}
\tilde{\rho}^{(N)} & = & \frac12 |S_{N-1}|  (\pi^2 t \tau)^{-\frac{N}{2}}   dx \, dy
\int_{S_{N-1}^+} dv D^{-\frac{N}{2}}
 e^{-\frac{1}{\tau D} (x^2 -2 v_N^2 xy + y^2(1+\tau t^{-1}))} \nonumber \\
&&   \hspace{.1in} \int_{R^{0|4(N-1)}}  d \phi d\bar{\phi} d\psi d\bar{\psi}
\left(v_N^2 + (y-x)(\hat{v}^T \bar{\phi}) (\hat{v}^T \phi) \right)
e^{-x\bar{\phi}^T\!\!\phi - y\bar{\psi}^T \!\! \psi} \nonumber \\
&&   \hspace{.2in}  e^{ - \frac{t}{2} (\bar{\phi}^T \bar{\psi}) (\phi^T \psi)  } \;
e^{ \frac{t^2}{2\tau D} (\bar{\phi}^T \bar{\psi}) (\hat{v}^T \phi) (\hat{v}^T \psi) } \;
  e^{ \frac{t}{\tau D}( x (\hat{v}^T \bar{\psi})(\hat{v}^T \psi)+ y (\hat{v}^T \bar{\phi})(\hat{v}^T \phi) )}.
\label{basepoint5}
\end{eqnarray}

\noindent
\textbf{Integration over Grassmann variables.}
%%%%%%%%%%%%%%%%%%%%%%%%%%%%%%%%%%%%%%%%%%%%%%%%%%
Let $R_v:\R^{N-1} \to \R^{N-1} $ be an orthogonal transformation that sends $\hat{v}$ to $(|\hat{v}|,0,\ldots,0)$.
We also change variables $\phi,\psi,\bar{\phi},\bar{\psi}$ to $R_v \phi, R_v \psi, R_v \bar{\phi}, R_v\bar{\psi}$.
Under this change $\hat{v}^T \phi$ becomes $|\hat{v}| \phi_1$ e.t.c. In terms of new variables the integral in
(\ref{basepoint5}) is
\begin{eqnarray}
\tilde{\rho}^{(N)} & = & \frac12 |S_{N-1}|  (\pi^2 t \tau)^{-\frac{N}{2}}
\int_{S_{N-1}^+} \! \! dv D^{-\frac{N}{2}}   e^{-\frac{1}{\tau D} (x^2 -2 v_N^2 xy + y^2(1+\tau t^{-1}))} \; I_5,
\label{basepoint6}
\end{eqnarray}
where
\begin{eqnarray}
I_5 &=&\int_{R^{0|4(N-1)}}  d \phi d \bar{\phi} d\psi d\bar{\psi}
\, \left(v_N^2 + (y-x)|\hat{v}|^2 \bar{\phi}_1 \phi_1\right) \nonumber\\
&& \hspace{.1in} e^{ -x\bar{\phi}^T\!\!\phi - y\bar{\psi}^T \!\! \psi} \;
e^{  - \frac{t}{2} (\bar{\phi}^T \bar{\psi}) (\phi^T \psi)} \;
e^{\frac{t^2 |\hat{v}|^2}{2\tau D}  (\bar{\phi}^T \bar{\psi}) \phi_1 \psi_1} \;
 e^{ \frac{t|\hat{v}|^2}{\tau D}( x \bar{\psi}_1 \psi_1 +  y \bar{\phi}_1 \phi_1)} \nonumber\\
 & = & v_{N}^2 \int_{R^{0|4(N-1)}}  d \phi d \bar{\phi} d\psi d\bar{\psi} \;
e^{\frac{(y-x)|\hat{v}|^2}{v_N^2} \bar{\phi}_1 \phi_1} \;
e^{ -x\bar{\phi}^T\!\!\phi - y\bar{\psi}^T \!\! \psi} \nonumber\\
&& \hspace{.1in}
e^{  - \frac{t}{2} (\bar{\phi}^T \bar{\psi}) (\phi^T \psi)} \;
e^{\frac{t^2 |\hat{v}|^2}{2\tau D}  (\bar{\phi}^T \bar{\psi}) \phi_1 \psi_1} \;
 e^{ \frac{t|\hat{v}|^2}{\tau D}( x \bar{\psi}_1 \psi_1 +  y \bar{\phi}_1 \phi_1)}.
\end{eqnarray}
Notice the appearance of fourth-order terms $ (\bar{\phi}^T \bar{\psi}) (\phi^T \psi) $ and
$  (\bar{\phi}^T \bar{\psi}) \phi_1 \psi_1$ in the above integral.
To deal with these we follow \cite{sommers2008general} and apply 
Hubbard-Stratonovich transformation
to convert the Berezin integral $I_5$ to an integral over both commuting 
and anti-commuting variables
but which is Gaussian with respect to the Grassmann variables, by using the identity
\begin{equation} \label{Strick}
e^{ ab } =  \int_{R^2} \frac{dz d\bar{z}}{\pi} \;
e^{ - z \bar{z}} e^{ a z} e^{ b \bar{z}}.
\end{equation}
In this identity, we have $z=x+iy, \, \bar{z} = x-iy $, and $dz d\bar{z} $ is shorthand for a Lebesgue integral over
$\R^2$ of this complex integrand.
The result is
\begin{eqnarray}
I_5 & = & v_{N}^2\int_{\R^2} \frac{dz d\bar{z}}{\pi}\int_{\R^2} \frac{dw d\bar{w}}{\pi}e^{-z\bar{z}-w\bar{w}}\int_{R^{0|4(N-1)}}
d \phi d \bar{\phi} d\psi d\bar{\psi}\nonumber\\
&& \hspace{.2in} e^{  (\frac{t}{2})^\alpha z (\bar{\phi}^T \bar{\psi}) - (\frac{t}{2})^\beta \bar{z} (\phi^T \psi)}
\; e^{(\frac{t^2 |\hat{v}|^2}{2\tau D})^\gamma  w (\bar{\phi}^T \bar{\psi}) + (\frac{t^2 |\hat{v}|^2}{2\tau D})^\delta
\bar{w} \phi_1 \psi_1}\nonumber\\
&& \hspace{.4in} e^{-x\bar{\phi}^T\!\!\phi - y\bar{\psi}^T \!\! \psi} \;
e^{ \frac{t|\hat{v}|^2}{\tau D}( x \bar{\psi}_1 \psi_1 +  y \bar{\phi}_1 \phi_1)}
\; e^{\frac{(y-x)|\hat{v}|^2}{v_N^2} \bar{\phi}_1 \phi_1},
\end{eqnarray}
where $\alpha+\beta=1$, $\gamma+\delta=1$. The final answer should not depend on $\alpha, \beta, \gamma, \delta$, which
provides us with a consistency check for the calculations below.
The internal Berezin integral breaks into the product of $N-1$ independent $4$-dimensional integrals, of which $(N-2)$
integrals are identical:
\bea
I_5 = v_{N^2}\int_{\R^2} \frac{dz d\bar{z}}{\pi}\int_{\R^2}
\frac{dw d\bar{w}}{\pi}e^{-z\bar{z}-w\bar{w}} \; (I_6)^{N-2} I_7
\eea
where
\bea
I_6 = \int_{R^{0|4}}  d \phi d \bar{\phi} d\psi d\bar{\psi} \;
e^{(\frac{t}{2})^\alpha z \bar{\phi} \, \bar{\psi}
- \left(\frac{t}{2}\right)^\beta \bar{z} \phi \psi} \;
e^{(\frac{t^2 |\hat{v}|^2}{2\tau D})^\gamma  w \bar{\phi} \, \bar{\psi}} \;
e^{-x\bar{\phi}\phi - y\bar{\psi}  \psi}
\eea
and
\begin{eqnarray}
I_7 &=&  \int_{R^{0|4}}  d \phi d \bar{\phi} d\psi d\bar{\psi} \;
e^{ (\frac{t}{2})^\alpha z \bar{\phi} \, \bar{\psi} - (\frac{t}{2})^\beta \bar{z} \phi \psi} \;
e^{(\frac{t^2 |\hat{v}|^2}{2\tau D})^\gamma  w \bar{\phi} \, \bar{\psi}
+(\frac{t^2 |\hat{v}|^2}{2\tau D})^\delta \bar{w} \phi \psi} \nonumber\\
& & \hspace{.4in} e^{ -x\bar{\phi} \, \phi - y \bar{\psi} \,  \psi} \;
e^{\frac{t|\hat{v}|^2}{\tau D}( x \bar{\psi}  \psi +  y \bar{\phi} \phi)} \;
e^{ \frac{(y-x)|\hat{v}|^2}{v_N^2} \bar{\phi} \phi}.
\end{eqnarray}
Each Berezin integral is Gaussian and can be evaluated using 
$ \int d \xi \exp(-\xi^T A \xi/2)=\mbox{Pf}(A) $.
Hence $I_6$ has value
\bea
\mbox{Pf}
\left( \begin{array}{cccc}
0 & x & - \left(\frac{t}{2}\right)^\beta \bar{z} & 0 \\
\cdot & 0 & 0 &  \left(\frac{t}{2}\right)^\alpha z+\left(\frac{t^2 |\hat{v}|^2}{2\tau D}\right)^{\gamma}w \\
\cdot & \cdot & 0 & y \\
\cdot & \cdot & \cdot & 0
\end{array} \right)=xy+\frac{t}{2}z\bar{z}+\left(\frac{t}{2}\right)^\beta \left(\frac{t^2 |\hat{v}|^2}{2\tau D}\right)^{\gamma}w\bar{z}\nonumber\\
\eea
and $I_7$ has value
\begin{eqnarray}
& & \hspace{-.5in} \mbox{Pf} \left( \begin{array}{cccc}
0 & \frac{x}{v_N^2}-y\frac{|\hat{v}|^2}{v_N^2}\frac{t+\tau}{\tau D} & - \left(\frac{t}{2}\right)^\beta \bar{z}+
\left(\frac{t^2 |\hat{v}|^2}{2\tau D}\right)^\delta \bar{w} & 0 \\
\cdot & 0 & 0 & \left(\frac{t}{2}\right)^\alpha z+\left(\frac{t^2 |\hat{v}|^2}{2\tau D}\right)^\gamma w \\
\cdot & \cdot & 0 & y-\frac{t|\hat{v}|^2}{\tau D}x \\
\cdot & \cdot & \cdot & 0
\end{array} \right)\nonumber\\
&=&\frac{1}{v_N^2}\left(x-y|\hat{v}|^2\frac{t+\tau}{\tau D}\right)\left(y-x|\hat{v}|^2 \frac{t}{\tau D}\right)+\frac{t}{2}z\bar{z}
-\frac{t^2|\hat{v}|^2}{2\tau D} w\bar{w}\nonumber\\
&& \hspace{.5in} -\left(\frac{t}{2}\right)^\alpha \left(\frac{t^2 |\hat{v}|^2}{2\tau D}\right)^\delta \bar{w} z
+\left(\frac{t}{2}\right)^\beta \left(\frac{t^2 |\hat{v}|^2}{2\tau D}\right)^\gamma w \bar{z}.
\end{eqnarray}

\noindent
\textbf{Integration over Hubbard Stratonovich variables.}
%%%%%%%%%%%%%%%%%%%%%%%%%%%%%%%%%%%%%%%%%%%%%%%%%%
To compute the integrals over complex variables $z,w$
notice that
\bea
\frac{dz d\bar{z}}{\pi}\int_C \frac{dw d\bar{w}}{\pi}e^{-z\bar{z}-w\bar{w}} z^k \bar{z}^l w^m \bar{w}^n
\eea
is zero for integers $k,l,m,n$ unless $k=l$ and $m=n$. Note that
$I_7$ contains terms in $\bar{w} z$ and $w \bar{z}$. So if
we expand $(I_6)^{N-2}$ as
\bea
%\left(xy+\frac{t}{2}z\bar{z}+\left(\frac{t}{2}\right)^\beta \left(\frac{t^2 |\hat{v}|^2}{2\tau D}\right)^{\gamma}w\bar{z}\right)^{N-2}&=&
\left(xy+\frac{t}{2}z\bar{z}\right)^{N-2} + (N-2)\left(xy+\frac{t}{2}z\bar{z}\right)^{N-3}
\left(\frac{t}{2}\right)^\beta \left(\frac{t^2 |\hat{v}|^2}{2\tau D}\right)^{\gamma}w\bar{z}+\ldots
\eea
we may ignore the other terms as they give zero contribution upon integration with respect to $w,z$.
Substituting the values of $I_6$ and $I_7$ we find
\begin{eqnarray}
&& I_5  =  v_{N}^2\int_C\frac{dz d\bar{z}}{\pi}\int_C \frac{dw d\bar{w}}{\pi}e^{-z\bar{z}-w\bar{w}}  \nonumber\\
&&
\left( \frac{1}{v_N^2}\left(x-y|\hat{v}|^2\frac{t+\tau}{\tau D}\right)\left(y-x|\hat{v}|^2
\frac{t}{\tau D}\right)+\frac{t}{2}z\bar{z}-\frac{t^2|\hat{v}|^2}{2\tau D} w\bar{w}\ \right) \left(xy+\frac{t}{2}z\bar{z}\right)^{N-2} \nonumber\\
&& -v_{N}^2(N-2)\frac{t^3|\hat{v}|^2}{4\tau D}\int_C\frac{dz d\bar{z}}{\pi}\int_C \frac{dw d\bar{w}}{\pi}e^{-z\bar{z}-w\bar{w}}
\left(xy+\frac{t}{2}z\bar{z}\right)^{N-3}w\bar{w}z\bar{z}.
\end{eqnarray}
These integrals can be computed using the formulae:
\begin{eqnarray}
\frac{1}{\pi}    \int_{R^2} dz d\bar{z}  e^{ - z \bar{z}} (\alpha +  |z|^2)^{N}
 =   \int^{\infty}_0 dr e^{-r}   (\alpha +  r)^{N} \nonumber \\
 = \sum_{k=0}^{N}{N \choose k} \alpha^k
\int^{\infty}_0 dr e^{-r}  r^{N-k}
 =   E_N(\alpha) \label{en1}
\end{eqnarray}
where $E_{N}(x)$ is the exponential polynomial of degree $N$ multiplied by $N!$,
\bea
E_N(x) = N! \, \sum_{k=0}^N \frac{x^k}{k!}.
\eea
Similarly
\begin{eqnarray}
\frac{1}{\pi}    \int_{R^2} dz d\bar{z} |z|^2 e^{ - z \bar{z}} (\alpha + |z|^2)^{N}
& = &  E_{N+1}(\alpha) - \alpha  E_N(\alpha) \nonumber \\
& = & (N+1) E_N(\alpha) - N \alpha E_{N-1}(\alpha). \label{en2}
\end{eqnarray}
So, integrating over $w,\bar{w}$ and pulling powers of $t$ out,
\begin{eqnarray}
I_5 & = & v_{N}^2\left(\frac{t}{2}\right)^{N-2}\int_C\frac{dz d\bar{z}}{\pi}e^{-z\bar{z}}\nonumber \\
&& \left( \frac{1}{v_N^2}\left(x-y|\hat{v}|^2\frac{t+\tau}{\tau D}\right)\left(y-x|\hat{v}|^2 \frac{t}{\tau D}\right)
+\frac{t}{2}z\bar{z}-\frac{t^2|\hat{v}|^2}{2\tau D} \ \right)  \left(\frac{2xy}{t}+z\bar{z}\right)^{N-2} \nonumber\\
&& -v_{N}^2(N-2)\left(\frac{t}{2}\right)^{N-3}\frac{t^3|\hat{v}|^2}{4\tau D}\int_C\frac{dz d\bar{z}}{\pi}e^{-z\bar{z}}
\left(\frac{2xy}{t}+z\bar{z}\right)^{N-3}z\bar{z}.
\end{eqnarray}
Finally, integrating over $z, \bar{z}$, we find, using (\ref{en1}) and (\ref{en2}),
\begin{eqnarray}
I_5 & = & \left(\frac{t}{2}\right)^{N-2}
\left(x-y|\hat{v}|^2\frac{t+\tau}{\tau D}\right)\left(y-x|\hat{v}|^2 \frac{t}{\tau D}\right)
E_{N-2}\left(\frac{2xy}{t}\right) \nonumber\\
& & \hspace{.1in}   -  \left(\frac{t}{2}\right)^{N-2} \frac{t^2|\hat{v}|^2v_N^2}{2\tau D}
E_{N-2}\left(\frac{2xy}{t}\right)\nonumber\\
&& \hspace{.2in} +v_{N}^2\left(\frac{t}{2}\right)^{N-1}\left((N-1) E_{N-2}\left(\frac{2xy}{t}\right)
- (N-2) \frac{2xy}{t} E_{N-2}\left(\frac{2xy}{t}\right)  \right)
\nonumber\\
&& \hspace{.3in} - v_{N}^2(N-2)\left(\frac{t}{2}\right)^{N-3}\frac{t^3|\hat{v}|^2}{4\tau D}
\left(E_{N-2}\left(\frac{2xy}{t}\right) - \frac{2xy}{t}E_{N-3}\left(\frac{2xy}{t}\right)  \right) \nonumber\\
  & = & \left(\frac{t}{2}\right)^{N-2}
\left(x-y|\hat{v}|^2\frac{t+\tau}{\tau D}\right)\left(y-x|\hat{v}|^2 \frac{t}{\tau D}\right)
E_{N-2}\left(\frac{2xy}{t}\right) \nonumber\\
& & \hspace{.1in}  + \frac{(N-1) v_N^2}{D} \left(\frac{t}{2}\right)^{N-1}  E_{N-2}\left(\frac{2xy}{t}\right) \nonumber\\
&& \hspace{.2in}
- \frac{(N-2) v_N^2}{D} \left(\frac{t}{2}\right)^{N-1} \frac{2xy}{t}E_{N-3}\left(\frac{2xy}{t}\right).
\end{eqnarray}
Substituting this into (\ref{basepoint6}), and using $\tau D = \tau + t | \hat{v}|^2$,
 gives us the following exact representation for the modified
density $\tilde{\rho}^{(N)}$ in the form of the integral over the $(N-1)$-dimensional sphere
\begin{eqnarray}
&& \hspace{-.4in} \tilde{\rho}^{(N)}_{t,t+\tau}(y,x) \nonumber \\
& = & \frac{2 |S_{N-1}|}{t^2 (2\pi)^N}
\int_{S_{N-1}^+} \! \! dv  \left(\tau t^{-1} + |\hat{v}|^2\right)^{-\frac{N}{2}}
e^{-\frac{1}{\tau + t |\hat{v}|^2} (x^2 -2 v_N^2 xy + y^2(1+\tau t^{-1}))} \nonumber \\
&& \left[
\left(x-y|\hat{v}|^2\frac{t+\tau}{\tau + t |\hat{v}|^2}\right)\left(y-x|\hat{v}|^2 \frac{t}{\tau + t |\hat{v}|^2}\right)
E_{N-2}\left(\frac{2xy}{t}\right) \nonumber \right.\\
& & \hspace{.1in} \left. + \frac{(N-1) t \tau v_N^2}{2(\tau+t|\hat{v}|^2)}  E_{N-2}\left(\frac{2xy}{t}\right)
- \frac{(N-2) t \tau v_N^2}{2(\tau + t |\hat{v}|^2)}  \frac{2xy}{t}E_{N-3}\left(\frac{2xy}{t}\right) \right].
\label{integral}
\end{eqnarray}
The integrand
depends only on $v_N$ and $|\hat{v}|^2 = 1-v_N^2$ and so the integral can be further reduced,
as we do in the following section, to a single integral over $v_N \in [0,1]$.
Moreover we will obtain very much simpler expressions for both $\tilde{\rho}$ and the spin-spin correlation
function valid in the limit $N\rightarrow \infty$.
%%%%%%%%%%%%%%%%%%%%%%%%%%%%%%%%%%%%%%%%%%%%%%%%%%%%%%%%
\section{Analysis of the large $N$ limit of the correlation function.} \label{largeN}
%%%%%%%%%%%%%%%%%%%%%%%%%%%%%%%%%%%%%%%%%%%%%%%%%%%%%%%%
A direct evaluation of (\ref{integral}) in the large-$N$ limit leads to the following conclusion:
\bea\label{eq87}
\lim_{N\rightarrow \infty}\tilde{\rho}^{(N)}_{t,t+\tau}(y,x)=0,
\eea
for any $\tau>0$.
%Via (\ref{inv2}) this implies that
%\bea
%\lim_{N\rightarrow \infty} \mathbb{E} \left(s_x(M_{t+\tau})s_y(M_{t})\right) = 0.
%\eea
%{\bf Theorem \ref{thm1} is proved.}
To study short-scale
correlations, we consider the scaling limit $\tau=\frac{T}{N}$ and $N \to \infty$.
We will calculate
\bea
\tilde{\rho}_{t,T}(y,x) = \lim_{N\rightarrow \infty}\tilde{\rho}^{(N)}_{t,t+T/N}(y,x).
\eea
It is possible to calculate $\tilde{\rho}_{t,T}(y,x)$ directly. The calculation is somewhat lengthy, but the answer
turns out to be, as expected, translational invariant. Assuming this translational invariance, it is possible to
calculate $\tilde{\rho}_{t,T}(y,x)$ as $\tilde{\rho}_{t,T}(0,x-y)$, which is a similar but less messy task, and
this is the calculation which we present here.

Firstly, when $y=0$ the expression (\ref{integral}) simplifies considerably.
Furthermore, using $v_N^2 + | \hat{v}|^2 = 1$, the integrand is
invariant under rotations in the equatorial plane and we may exploit the identity
\bea
\int_{S_{N-1}^+} dv \, g(v_N) = |S_{N-2}|  \int^1_{0} dz \, g(z) (1-z^2)^{\frac{N-3}{2}}.
\eea
We find
\begin{eqnarray}
\tilde{\rho}^{(N)} _{t,t+T/N}(0,x) &=&  c_1(N) \int^1_0 dz \, (1-z^2)^{-5/2}
e^{- \frac{x^2}{(T/N)+t(1-z^2)}} \\
& & \hspace{.2in}  \left(1+ \frac{T}{Nt(1-z^2)} \right)^{-\frac{N}{2}-1}
 \left( z^2 - \frac{2N x^2 (1-z^2)}{(N-1)T} \right)\nonumber
\end{eqnarray}
where
\bea
c_1(N) = \frac{(N-1)!}{(2 \pi)^N} |S_{N-1}| |S_{N-2}| \frac{T}{Nt^2}.
\eea
Using the volume formula $|S_{N-1}|=\frac{N \pi^{\frac{N}{2}}}{\Gamma(\frac{N}{2}+1)}$, and
asymptotic $\Gamma(z) \sim z^z z^{-1/2} e^{-z}$ for large $z$, it is straightforward to check that
$c_1(N) \to T/(\pi t^2)$ as $N \to \infty$. Moreover one may justify the
passage of the limit inside the $dz$ integral, for example by dominated convergence using
\bea
\left(1+ \frac{T}{Nt(1-z^2)} \right)^{-\frac{N}{2}} \leq (1-z^2)^2 \frac{8t^2}{T^2}.
\eea
This limit gives
\bea
\tilde{\rho}_{t,T} (0,x) = \frac{T}{\pi t^2}  \int^1_0 dz \, (1-z^2)^{-5/2}
e^{- \frac{x^2}{t(1-z^2)}}  \, e^{-\frac{T}{2t(1-z^2)}}
 \left( z^2 - \frac{2x^2 (1-z^2)}{T} \right).
\eea
The substitution $(1-z^2) = (1+w)^{-1} $ and $dz = \frac12 w^{-1/2} (1+w)^{-3/2} dw$ yields
\begin{eqnarray}
&& \hspace{-.4in} \tilde{\rho}_{t,T} (0,x) \nonumber \\
& = & \frac{T}{2\pi t^2} e^{- \left(\frac{x^2}{t}+\frac{T}{2t}\right)}
\int_{0}^{\infty} dw \, e^{-\left(\frac{x^2}{t}+\frac{T}{2t}\right) w}
\left(w^{1/2}-2\frac{x^2}{T}w^{-1/2}\right) \nonumber \\
& = & \frac{T}{2 \sqrt{\pi} t^2}
 e^{- \left(\frac{x^2}{t}+\frac{T}{2t}\right)}
\left(\frac12 \left(\frac{x^2}{t}+\frac{T}{2t}\right)^{-3/2}
- \frac{2x^2}{T}\left(\frac{x^2}{t}+\frac{T}{2t}\right)^{-1/2}\right)  \label{answer1}
\end{eqnarray}
which is the final answer for the scaling limit of the modified two-point density
in the large $N$ scaling limit.
%Note the exponential decay at time scales of order $O(1)$,
%which is consistent with the de-correlation time for $gl(N)$-Brownian motions
%being of order of $1/N$.
Note that the limit
\bea
\lim_{T \downarrow 0} \tilde{\rho}_{t,T}(0,x)  =-\frac{1}{\sqrt{\pi}}\frac{|x|}{t^{3/2}} e^{-\frac{x^2}{t}}
+ \frac{1}{\sqrt{\pi t}} \delta(x)
\eea
is consistent with the one-dimensional correlation function  obtained in \cite{borodin2009ginibre},
\cite{forrester2007eigenvalue} (the
presence of the delta function at zero is expected, see \cite{ECP2133} for the explanation). Let
\bea
R(t,y;T,x) = \lim_{N\rightarrow \infty} \mathbb{E} \left(s_{y}(M_t)s_{x}(M_{t+\frac{T}{N}})\right)
\eea
be the scaling limit of the unmodified two-time spin-spin correlation function.  Using (\ref{inv2}), this is given by
\begin{eqnarray}
R(t,y;T,x+y) 
& = & - \lim_{N \to \infty} 4 \int_{y}^{\infty} \! dy' \int_{-\infty}^{x+y} \! dx' \,  \tilde{\rho}^{(N)}_{t,t+T/N}(y',x') \nonumber \\
& = & - 4 \int_{y}^{\infty} \! dy' \int_{-\infty}^{x+y} \! dx' \,  \tilde{\rho}_{t,T}(y',x') \nonumber \\
& = &  -4 \int_{y}^{\infty} \! dy' \int_{-\infty}^{x+y} \! dx' \,  \tilde{\rho}_{t,T}(0,x'-y') \nonumber \\
& = & 4 \int_{-\infty}^{x} dz (z-x) \tilde{\rho}_{t,T}(0,z).\label{care}
\end{eqnarray}
Justifying the passage of the limit $N \to \infty$ inside the integral in (\ref{care}) needs a little care,
and we delay the argument until the end of this section.
Substituting (\ref{answer1}) into the above, we find that
\bea
R(t,y;T,x+y) = \frac{4}{\sqrt{\pi} t} \int_{x}^\infty dz (x-z)e^{-\alpha(z)} \left(\frac{1}{2}r\alpha^{-3/2}(z)
-(\alpha(z)-r)\alpha^{-1/2}(z)\right),\nonumber\\
\eea
where $r=\frac{T}{2t}$, $\alpha(z)=\frac{z^2}{t}+r$. Since $\tilde{\rho}_{t,T}(x,0)$ is even in
$x$ so too is $R(t,y;T,x+y)$. Taking $x>0$,
the change of variables $w = \alpha(z) $, so that $z/t^{1/2} = (\alpha(z)-r)^{1/2}$, gives
\begin{eqnarray}
&& \hspace{-.3in} R(t,y;T,x+y) \nonumber\\
& = & - \frac{2}{\sqrt{\pi}} \int^{\infty}_{\alpha(x)} dw \, e^{-w}
\left( \frac12 r w^{-3/2} - (w-r) w^{-{1/2}} \right) \nonumber\\
& & \hspace{.2in} + \frac{2x}{\sqrt{\pi t}}  \int^{\infty}_{\alpha(x)} dw \, e^{-w}
\left( \frac12 r (w-r)^{-1/2} w^{-3/2} - (w-r)^{1/2} w^{-{1/2}} \right)  \nonumber\\
& = & - \frac{2r}{\sqrt{\pi}} e^{-\alpha(x)} \alpha(x)^{-1/2} +
 \frac{2}{\sqrt{\pi}} \int^{\infty}_{\alpha(x)} dw \, e^{-w} w^{1/2} \nonumber\\
& & \hspace{.2in} + \frac{2x}{\sqrt{\pi t}}  \int^{\infty}_{\alpha(x)} dw \,
\frac{d}{dw} \left( e^{-w}(w-r)^{1/2} w^{-{1/2}} \right)  \nonumber\\
& = &  \frac{1}{\sqrt{\pi}}  \int^{\infty}_{\alpha(x)} dw \, e^{-w} w^{-1/2},
\end{eqnarray}
where the last two equalities follow by integration by parts, and all the
boundary terms cancel.
This final integral is the same as the complementary error function giving
\bea
R(t,y;T,x+y)=\mbox{\erfc}\left(\sqrt{\frac{x^2}{t}+\frac{T}{2t}}\right).
\label{answer2}
\eea
{\bf Theorem \ref{thm2} is proved.}

We now return to complete the justification of (\ref{care}).
We need to pass to the limit to show
\begin{equation} \label{form1}
\int^{\infty}_y dy' \int^x_{-\infty} dx' \tilde{\rho}^{(N)}_{t,t+T/N}(y',x') \to
\int^{\infty}_y dy' \int^x_{-\infty} dx' \tilde{\rho}_{t,T}(y',x').
\end{equation}
By $y \to -y, x \to -x$ symmetry we may suppose $y \geq 0$.
Using the fact that $\tilde{\rho}^{(N)}$ is a mixed derivative we know
\begin{equation} \label{symtrick}
\int_{\R} \tilde{\rho}^{(N)}_{t,\tau}(x,y) dx = \int_{\R} \tilde{\rho}^{(N)}_{t,\tau}(x,y) dy = 0.
\end{equation}
The same property holds for the limit $\tilde{\rho}_{t,T}(x,y)$. When $y \geq 0$ and $x \geq 0$ it is
convenient to use (\ref{symtrick}) to rewrite (\ref{form1}) as
\begin{eqnarray}
&& \hspace{-.3in}
\int^{\infty}_y dy' \int^0_{-\infty} dx' \tilde{\rho}^{(N)}_{t,t+T/N}(y',x')
- \int_{-\infty}^y dy' \int^x_{0} dx' \tilde{\rho}^{(N)}_{t,t+T/N}(y',x')
\nonumber\\
& \to &
\int^{\infty}_y dy' \int^0_{-\infty} dx' \tilde{\rho}_{t,T}(y',x')
- \int_{-\infty}^y dy' \int^x_{0} dx \tilde{\rho}_{t,T}(y',x').
\label{form2}
\end{eqnarray}
This ensures that the region of integration has only a bounded region
where $xy >0$, namely $[0,x] \times [0,y]$.

We will replace $\tilde{\rho}^{(N)}$
by a further modification $\hat{\rho}^{(N)}$, where we replace each occurrence of
the truncated exponentials $E_N(z)$ by the un-truncated
exponential $N! e^z$. Thus we define (compare with (\ref{integral}))
\begin{eqnarray}
&& \hspace{-.4in} \hat{\rho}^{(N)}_{t,t+\tau}(y,x) \nonumber \\
& = & \frac{2 |S_{N-1}|}{t^2 (2\pi)^N}
\int_{S_{N-1}^+} \! \! dv  \left(\tau t^{-1} + |\hat{v}|^2\right)^{-\frac{N}{2}}
e^{-\frac{1}{\tau + t |\hat{v}|^2} (x^2 -2 v_N^2 xy + y^2(1+\tau t^{-1}))} e^{\frac{2xy}{t}} \nonumber \\
&& \left[
\left(x-y|\hat{v}|^2\frac{t+\tau}{\tau + t |\hat{v}|^2}\right)\left(y-x|\hat{v}|^2 \frac{t}{\tau + t |\hat{v}|^2}\right)
(N-2)!  \nonumber \right.\\
& & \hspace{1.3in} \left. + \frac{(N-1)! \, t \tau v_N^2}{2(\tau+t|\hat{v}|^2)}
- \frac{(N-2)! \, t \tau v_N^2}{2(\tau + t |\hat{v}|^2)}  \frac{2xy}{t}  \right].
\label{integral2}
\end{eqnarray}
This does not change the point-wise limit, and we still have
\bea
\tilde{\rho}_{t,T}(y,x) = \lim_{N\rightarrow \infty} \hat{\rho}^{(N)}_{t,t+T/N}(y,x).
\eea
We need a uniform bound for $\hat{\rho}^{(N)}_{t,t+T/N}(y,x) $. We write
$C(t,T,\ldots)$ for a constant whose value may change but whose dependency is indicated. The key term
is the exponential: when $xy \leq 0$ we may bound
\bea
e^{-\frac{1}{TN^{-1} + t |\hat{v}|^2} \left(x^2 -2 v_N^2 xy + y^2(1+T N^{-1}t^{-1})\right)}
 e^{\frac{2xy}{t}}
\leq e^{-\frac{x^2+y^2}{T+t}};
\eea
and when $xy \geq 0$ we use
\begin{eqnarray}
e^{-\frac{1}{TN^{-1} + t |\hat{v}|^2} \left(x^2 -2 v_N^2 xy + y^2(1+T N^{-1}t^{-1})\right)}
 e^{\frac{2xy}{t}}
& = & e^{-\frac{1}{TN^{-1} + t |\hat{v}|^2} \left((x-y)^2 + (y^2-2xy)T N^{-1}t^{-1}\right)}\nonumber \\
& \leq & e^{-\frac{(x-y)^2}{T+t}} e^{\frac{2xy}{t}}.
\end{eqnarray}
The terms in the square brackets in the integrand in  $\hat{\rho}^{(N)}_{t,t+T/N}(y,x) $ can be simply bounded
by $(N-2)! C(t,T) (T+x^2+y^2)|\hat{v}|^{-2}$, and carrying out the integral over $S_{N-1}^+$ yields the bound
\begin{eqnarray} \label{hatbound}
\hat{\rho}^{(N)}_{t,t+T/N}(y,x)
\leq
C(t,T) (T+x^2+y^2)
\left\{
\begin{array}{ll}
e^{-\frac{x^2+y^2}{T+t}} & \mbox{if $xy \leq 0$,} \\
e^{-\frac{(x-y)^2}{T+t}} e^{\frac{2xy}{t}} & \mbox{if $xy \geq 0$.}
\end{array} \right.
\end{eqnarray}
Notice the disappearance of the $N$-dependence from the above estimate.
When $y \geq 0$ and $x \leq 0$ the limit
\bea
\int^{\infty}_y dy' \int^x_{-\infty} dx' \hat{\rho}^{(N)}_{t,t+T/N}(y',x') \to
\int^{\infty}_y dy' \int^x_{-\infty} dx' \tilde{\rho}_{t,T}(y',x')
\eea
follows from (\ref{hatbound}) by dominated convergence. When $y \geq 0$ and $x \geq 0$ we switch
to the form (\ref{form2})
and obtain, also by dominated convergence using (\ref{hatbound}),
\begin{eqnarray}
&& \hspace{-.3in}
\int^{\infty}_y dy' \int^0_{-\infty} dx' \hat{\rho}^{(N)}_{t,t+T/N}(y',x')
- \int_{-\infty}^y dy' \int^x_{0} dx' \hat{\rho}^{(N)}_{t,t+T/N}(y',x')
\nonumber\\
& \to &
\int^{\infty}_y dy' \int^0_{-\infty} dx' \tilde{\rho}_{t,T}(y',x')
- \int_{-\infty}^y dy' \int^x_{0} dx \tilde{\rho}_{t,T}(y',x').
\end{eqnarray}
It remains to show the error from approximating $\tilde{\rho}^{(N)}$ by $\hat{\rho}^{(N)}$
is negligible. However this is a simpler task, using the fact that $E_N(z)/N!$ is close to $e^z$
on a ball of radius $N^{1/2}$, and using the simple bound $|E_N(z)| \leq N! e^{|z|}$ outside this ball.
We omit the details.

Theorem \ref{thm1} follows from (\ref{eq87}) using the integration formula (\ref{inv2}). The interchange
of $N\rightarrow \infty$ limit and integration can be justified using dominated convergence argument
similar to the one given above. The conclusion is
\bea
\lim_{N\rightarrow \infty} \mathbb{E} \left(s_x(M_{t+\tau})s_y(M_{t})\right) = 0.
\eea
{\bf Theorem \ref{thm1} is proved.}

%
%%%%%%%%%%%%%%%%%%%%%%%%%%%%%%%%%%%%%%%%%%%%%%
%%%%%%%%%%%%%%%%%%%%%%%%%%%%%%%%%%%%%%%%%%%%%%

\section{Fixed time multi-point density functions of real eigenvalues via spin variables.}
\label{oned}
As another application of the formalism of spin variables
we will show how to
re-derive the result of Borodin-Sinclair-Forrester-Nagao \cite{borodin2009ginibre}, \cite{forrester2007eigenvalue} for all the
one-dimensional densities of real eigenvalues of the real Ginibre ensemble in the limit $N\rightarrow \infty$.
And moreover it prepares grounds for the future study of multi-time correlation functions.

Let us fix time $t=1$ and a set of  multiple
space points $x_1<x_2<\ldots<x_K$. The answer for an arbitrary time $t$ can be obtained from the $t=1$ answer by the
diffusive re-scaling $x\rightarrow \frac{x}{\sqrt{t}}$. As we have done before, we will first compute a modified density
$ \tilde{\rho}^{(N)} =  \tilde{\rho}^{(N)} (x_1,x_2,\ldots,x_K) $ defined by
\bea
\tilde{\rho}^{(N)} (x_1,x_2,\ldots,x_K) \, dx_1 \ldots dx_K =
\mathbb{E}_N \left( \prod_{k=1}^K s_{x_k}(M_1) \Lambda^{M_1}(dx_k)\right)
\eea
where
where the subscript $N$ on the $\mathbb{E}_{N}$ means that the averaging occurs over
the $N \times N$ Ginibre distribution.
As before, we choose right hand limits for the density, that is  where the intervals $dx_k$ denote
infinitesimal intervals just to the right of the point $x_k$.

The answers are zero for odd $K$, so we take {\bf even $K$ throughout}. 
Moreover it is convenient to consider only
{\bf even $N$ throughout} (which avoids us tracking various $\pm$ signs).

\noindent
\textbf{Equivalence with a correlation function of characteristic polynomials.}
%%%%%%%%%%%%%%%%%%%%%%%%%%%%%%%%%%%%%%%%%%%%%%%%%%
The integral over the Gaussian density
\bea
\tilde{\rho}^{(N)} (x_1,x_2,\ldots,x_K) \, dx_1 \ldots dx_K = \int_{R^{N^2}} dM \gamma_1(M)
\prod_{k=1}^K s_{x_k}(M)  \Lambda^{M}(dx_k)
\eea
can be treated using the Edelman transform (\ref{edt}, \ref{edt1}, \ref{pvform}) for the eigenvalue lying in $dx_K$ as in
section \ref{sec:dervtns}, and after integrating over the half sphere $S_{N-1}^+$, we obtain
\begin{eqnarray}
 && \hspace{-.3in} \tilde{\rho}^{(N)}  (x_1,x_2,\ldots,x_K)  \, dx_1 \ldots dx_{K-1} \\
   & = &   \frac12 |S_{N-1}| \pi^{-\frac{N-1}{2}} e^{-x_K^2}
  \mathbb{E}_{N-1} \left(\det\left(M_1-x_K I \right) \prod_{k=1}^{K-1} s_{x_k}(M_1) \Lambda^{M_1}(dx_k)\right).\nonumber
\end{eqnarray}
Another Edelman transform about the eigenvalue lying in $dx_{K-1}$ yields
\begin{eqnarray}
 && \hspace{-.3in} \tilde{\rho}^{(N)}  (x_1,x_2,\ldots,x_K)  \, dx_1 \ldots dx_{K-2} \nonumber\\
& = &  \frac14 |S_{N-1}| |S_{N-2}|
\pi^{-\frac{N-1}{2}- \frac{N-2}{2}} e^{-x^2_K-x^2_{K-1}} (x_{K-1}-x_{K}) \\
&&  \hspace{.2in} \mathbb{E}_{N-2} \left(\det\left(M_1-x_K I \right) \det\left(M_1-x_{K-1} I \right)
\prod_{k=1}^{K-2} s_{x_k}(M_1) \Lambda^{M_1}(dx_k)\right).\nonumber
\end{eqnarray}
A further $(K-2)$ applications of Edelman transform will lead
to the following expression for the modified density:
\begin{eqnarray}
&& \hspace{-.3in} \tilde{\rho}^{(N)}(x_1,x_2,\ldots, x_K )  \nonumber \\
& = & \frac{\Delta(\mathbf{x}) }{2^K} \prod_{k=1}^{K} \left( |S_{N-k}|
\pi^{-\frac{N-k}{2}} e^{-x_k^2} \right)  \mathbb{E}_{N-K} \! \left(\prod_{m=1}^K\det\left(M_1-x_m I \right)\right)
\label{app:md}
\end{eqnarray}
where $\Delta (\mathbf{x})=\prod_{1\leq i < j \leq K} (x_j-x_i)$ is the Vandermonde determinant.
Therefore, the problem of computing the modified density has been reduced to the computation of the expectation
of the product of characteristic polynomials of the random matrix $M_1$.

Averaging products of characteristic polynomials is a well-studied problem in random matrix theory, see \cite{handbook}
for a review. In principle, we could have stopped here by pointing out that the desired Pfaffian expression for the correlation
functions of real eigenvalues follow for example from (\ref{app:md}) and the Pfaffian formulae of \cite{guhr_kieburg}.

The reason for pressing on with the calculation is three-fold: firstly, the method presented below allows for a straightforward generalization to the multi-time case which will report on in \cite{tribe_fc}. Secondly, our calculation 
uncovers an interesting new integral representation for (\ref{app:md}) in terms of a seemingly new exactly localizing integral over $U(K)/USp(K)$, see Theorem \ref{thm3}. This integral plays an important role in the analysis of multi-time correlation
functions for the Ginibre evolution and we felt that the community should know about it. Thirdly, our aim is to present a rigorous
route from spin-spin correlation functions to Borodin-Sinclair-Forrester-Nagao result. Thus we are forced to spend some time
on a, perhaps dull, analysis of convergence of (\ref{app:md}) in the large-$N$ limit in order to prove 
that corrections to the leading term (given
by (\ref{temp11} below) do indeed vanish as $N\rightarrow \infty$.

\noindent
\textbf{Integral representation for product of characteristic polynomials.}
%%%%%%%%%%%%%%%%%%%%%%%%%%%%%%%%%%%%%%%%%%%%%%%%%%
Such a computation has been carried out in \cite{sommers2008general}
by exploiting Berezin integrals. As it is rather close to the calculation
already detailed in section \ref{sec:dervtns} we will just present the final answer.
\bea
&& \hspace{-.3in} \mathbb{E}_{N} \left(\prod_{m=1}^K \det \left(M_1-x_m I\right)\right) \nonumber \\
& = &  \prod_{1\leq p< q \leq K} \left[ \int_{\mathbf{\R^2}} \frac{dz_{pq} d\bar{z}_{pq}}{\pi}
e^{-|z_{pq}|^2} \right] Pf
\left( \begin{array}{cc}
\frac{1}{\sqrt{2}}Z & X  \\
-X & \frac{1}{\sqrt{2}}Z^{\dagger}
\end{array} \right)^N. \label{app:int1}
\eea
Here each $dz_{pq} d\bar{z}_{pq}$ is shorthand for Lebesgue measure on $\R^2$ and arises
from repeated use of the Hubbard Stratonovich transform; the matrix $X$ is a diagonal
$K\times K$ matrix with entries $(x_1, x_2, \ldots x_K)$; and $Z$ is a skew symmetric
complex $K\times K$ matrix:
\begin{eqnarray}
Z_{ij}=\left\{ \begin{array}{cc}
z_{ij} & i>j,  \\
0 &  i=j,\\
-z_{ij}& i<j.
\end{array} \right.
\end{eqnarray}
The right hand side of expression (\ref{app:int1}) can be neatly written as a matrix integral:
\bea\label{app:int2}
 \pi^{-\frac{K(K-1)}{2}}
 \int_{Q^{(K)}} \lambda(dZ,dZ^{\dagger}) e^{-\frac{1}{2} Tr ZZ^{\dagger}}
Pf
\left( \begin{array}{cc}
\frac{1}{\sqrt{2}}Z & X  \\
-X & \frac{1}{\sqrt{2}}Z^{\dagger}
\end{array} \right)^N,
\eea
where $Q^{(K)}=\{ Z\in \mathbf{C}^{K\times K}\mid Z^T=-Z^T\}$ is the space of skew-symmetric
complex matrices, $\lambda (Z,Z^{\dagger})$ is the Lebesgue measure on $Q^{(K)}$ as described above.

It is worth noting here, that the application  of the standard Hubbard-Stratonovich transformation
leads to convergent integrals over commuting variables. This is due to the fact that it is applied to exponentials
depending on anti-commuting variables. For the multi-time case this doesn't happen and one has to use
the so called super-bosonization technique instead, see e.g. \cite{mz1}, \cite{mz2}.

Note that the dimension of the integral
in the right hand side of (\ref{app:int2}) is $N$-independent. The size of the original
matrix only enters the integral as the power of the Pfaffian in the integrand. This
allows one to calculate the large $N$-limit of (\ref{app:int2}) using the Laplace method.
To facilitate the application of asymptotic methods, we re-scale the integration variables
using $(Z,\dZ) \rightarrow \sqrt{N} (Z,\dZ)$, which gives
\begin{equation}
\label{app:int3}
\mathbb{E}_{N} \left(\prod_{m=1}^K\det\left(M_1-x_m\right)\right)
 =  \pi^{-\frac{K(K-1)}{2}} 2^{-\frac{NK}{2} }N^{\frac{NK}{2}} N^{\frac{K(K-1)}{2}} J_N
\end{equation}
where
\bea
J_N = \int_{Q^{(K)}} \lambda(dZ,dZ^{\dagger}) e^{-\frac{N}{2} Tr ZZ^{\dagger}}
Pf
\left( \begin{array}{cc}
Z & \sqrt{\frac{2}{N}}X  \\
-\sqrt{\frac{2}{N}}X & Z^{\dagger}
\end{array} \right)^N \!\!\! \!\!\!.
\eea
The integrand in $J_N$ is now of the form $\exp (NF_N(Z))$, where $F_N$ is a slow function of $N$. (In the sense
that $F$ and its derivatives converge in the limit $N\rightarrow \infty$.)

\noindent
\textbf{Applying the Laplace method.}
%%%%%%%%%%%%%%%%%%%%%%%%%%%%%%%%%%%%%%%%%%%%%%%%%%
We will show that the integral $J_N$ localizes onto the subset
\begin{equation} \label{CK}
C^{(K)} = \{Z \in Q^{(K)}\mid Z\dZ=I \}.
\end{equation}
To do this it is convenient to split $J$ into two parts,
\bea
J_{N,0} = \int_{Q^{(K)} \cap S} \qquad J_{N,1} = \int_{Q^{(K)} \setminus S}
\eea
where
\bea
S = \left\{ Z \mid  \mu_k(Z) \in \left[ 1/2, 2\right] \; \mbox{for $k=1,\ldots,K$} \right\}
\eea
and $(\mu_k(Z): 1 \leq k \leq K\}$ are the singular values of $Z$. Note that $S$ is compact and contains
only non-singular matrices.
We first bound the integral $J_{N,1}$, aiming to show that it is of smaller order
that $J_{N,0}$.
 For $c \geq 0$ write $(\lambda_k(c): 1 \leq k \leq 2K)$ for the singular values
of the matrix
\bea
\left( \begin{array}{cc}
Z & c X  \\
-cX & Z^{\dagger}
\end{array} \right)
\eea
One may bound the difference of two sets of singular values $(\mu_k(A))$ and $(\mu_k(B))$
via the operator norm bound $|\mu_k(A) - \mu_k(B)| \leq \|A-B\|$. Thus
\bea
\left| \lambda_k(c) - \lambda_k(0) \right| \leq \| cX\| = c x_*
\eea
where $x_* = \max_k |x_k|$. Then
\begin{eqnarray}
\left| Pf
\left( \begin{array}{cc}
Z & c X  \nonumber\\
- c X & Z^{\dagger}
\end{array}  \right) \right|
& = &
\left| \det
\left( \begin{array}{cc}
Z & c X  \nonumber\\
-c X & Z^{\dagger}
\end{array}\right)  \right|^{1/2} \nonumber\\
& = & \left| \prod_{k=1}^{2K} \lambda_k(c) \right|^{1/2} \nonumber\\
& \leq & \left|  \prod_{k=1}^{2K} \left( \lambda_k(0) + c x_*\right) \right|^{1/2} \nonumber\\
& = & \prod_{k=1}^K (\mu_k(Z) + cx_*).
\end{eqnarray}
Using this in $J_{N,1}$  we find
\begin{eqnarray}
J_{N,1} & \leq &  \int_{Q^{(K)}\setminus S} \lambda(dZ,dZ^{\dagger}) e^{-\frac{N}{2} Tr ZZ^{\dagger}}
\left| \prod_{k=1}^K \left(\mu_k(Z) + (2/N)^{1/2} x_* \right) \right|^N \nonumber\\
%& = &  \int_{Q^{(K)}\setminus S} \lambda(dZ,dZ^{\dagger})
%\exp - N \sum_{k=1}^K \left( \frac12 \mu_k^2(Z) - \ln(\mu_k(Z) +  (2/N)^{1/2}x_*)\right) \\
& = & \int_{Q^{(K)}\setminus S} \lambda(dZ,dZ^{\dagger})
e^{- Tr ZZ^{\dagger}}
e^{- N \sum_{k=1}^K H_N(\mu_k(Z))}
\end{eqnarray}
where
\bea
H_N(z) = \left(\frac12 - \frac1N \right) z^2 - \ln \left(z +   (2/N)^{1/2}  x_* \right).
\eea
Note that $H_N(z) \to H(z) = \frac12 z^2 - \ln(z)$ and that $H(z)$ has a minimal value of $H(1)=\frac12$.
On $Q^{(K)} \setminus S$ there must exist at least one singular value $\mu$ lying outside $[\frac12,2]$,
and for this value $H(\mu) \geq \frac12 + 2 \delta$ for an easily calculated $\delta >0$. Thus for large $N$, when
$Z \in Q^{(K)} \setminus S$, we have
\bea
\sum_{k=1}^K H_N(\mu_k(Z)) \geq \frac{K}{2} + \delta
\eea
and the value of $J_{N,1}$ is bounded by $C(K) e^{-NK/2} e^{-N \delta}$. This is
is exponentially smaller than that of $J_{N,0}$, which we will see is, to leading exponential order, $O(e^{-NK/2})$.

Next, we will calculate the asymptotic expansion of $J_{N,0}$ for large $N$.
The $N$th power of the Pfaffian in the integrand of $J_{N,0}$ can be simplified using the Taylor
expansion for the Pfaffian:
\bea\label{app:pte}
\frac{Pf\left( A+\frac{1}{\sqrt{N}}B \right)}{Pf\left(A\right)}&=&
1+\frac{1}{2\sqrt{N}}TrBA^{-1}\\
&+&\frac{1}{8N}\left(TrBA^{-1}TrBA^{-1}-2 Tr BA^{-1}BA^{-1}\right)+O\left(N^{-\frac{3}{2}}\right)\nonumber
\eea
While  we cannot pinpoint the exact reference for the original derivation of the above expansion, it can be easily derived
using the Berezin integral representation of the Pfaffian. It is interesting to note that unlike the analogous determinant
expansion formula, the series (\ref{app:pte}) contains finitely many terms. Using (\ref{app:pte}),
and the fact that the terms with $Tr(BA^{-1})$ are zero in our case,
we can re-write the $N$th power of the Pfaffian in the integrand of $J_{N,0}$ as follows:
\begin{eqnarray}
Pf^N
\left( \begin{array}{cc}
Z & \sqrt{\frac{2}{N}}X  \nonumber\\
-\sqrt{\frac{2}{N}}X & Z^{\dagger}
\end{array} \right)
&=& \det(Z\dZ)^{\frac{N}{2}}  \left(1+\frac{1}{N}Tr\left(Z^{\dagger}X Z X\right)+O(N^{-2})\right)^N\\
&=& \det(Z\dZ)^{\frac{N}{2}}  e^{Tr\left(Z^{\dagger}X Z X\right)}\left(1+O(N^{-1})\right)
\end{eqnarray}
This allows us to express $J_{N,0}$ in a form well suited for the application of Laplace formula:
\begin{equation}
\label{app:io}
J_{N,0} = \int_{Q^{(K)} \cap S} \Lambda(dZ,dZ^{\dagger}) e^{-\frac{N}{2} \left(Tr ZZ^{\dagger}-\ln \det\left(Z\dZ\right)\right)}
e^{Tr\left(Z^{\dagger}X Z X\right)} (1+O(N^{-1})).
\end{equation}
The fact that $S$ does not contain degenerate matrices, and the compactness of $S$ allows one
to pass the correction term $O(N^{-1})$ through the integral.
In the limit $N\rightarrow \infty$, the main contribution to (\ref{app:io}) comes from the neighborhood of the
points of global minimum of the function
\bea
F(Z)=Tr ZZ^{\dagger}-\ln \det\left(Z\dZ\right) = \sum_{k=1}^{K} \left( \mu_k^2(Z) - 2 \ln( \mu_k(Z) ) \right).
\eea
The global minimum value of $F$ is $K$ and it is attained on the set $C^{(K)}$ in (\ref{CK}),
namely the skew-symmetric unitary $K \times K$ matrices, which is
a smooth sub-manifold of $Q^{(K)}$.  We will show that
$C^{(K)}$ is a non-degenerate critical set, which means that the Hessian of $F$ has the maximal possible
rank at every point of $C^{(K)}$. Therefore we can use Laplace theorem \cite{erdelyi2010asymptotic} to
calculate the asymptotic expansion of $J_{N,0}$: let $(w,y)$ be local co-ordinates on $Q^{(K)}$ such that the sub-manifold $C^{(K)}$
is locally determined by the set of equations $y=0$; then
\begin{eqnarray}
&& \hspace{-.3in} \int_{Q^{(K)} \cap S} \Lambda(dZ,dZ^{\dagger}) e^{-\frac{N}{2} \left(Tr ZZ^{\dagger}-\ln \det\left(Z\dZ\right)\right)}
e^{Tr\left(Z^{\dagger}X Z X\right)} \nonumber \\
&=&e^{-NF \mid_{C^{(K)}}}\left(\frac{1}{\sqrt{2\pi N}}\right)^{\mbox{dim}(Q^{(K)})-\mbox{dim}(C^{(K)})}
\int_{C^{(K)}}\mu(dw) e^{Tr\left(Z^{\dagger}X Z X\right)}. \label{app:laplace}
\end{eqnarray}
Here  $\mu(dw)$ is the measure
on $C^{(K)}$ generated by the embedding $C^{(K)}\subset Q^{(K)}$ and integration over transverse co-ordinates $y$. Explicitly,
\begin{equation} \label{app:meas}
d\mu(w)=\frac{\rho(w,y=0)}{\sqrt{\det Hess(F)\mid_C^{(K)} (w)}}\prod_{k=1}^{dim(C^{(K)})}dw_k,
\end{equation}
where $\rho(w,y)$ is the density of Lebesgue
 measure $\Lambda(dZ, d\dZ)$ with respect to local co-ordinates $(w,y)$, the Hessian is defined as the matrix of second derivatives
 with respect to transverse co-ordinates $y$. In writing (\ref{app:laplace}) we used the fact that the critical manifold lies
a positive distance away from the boundary of $Q^{(K)} \cap S$. Noting that $F$ takes the value $K$ on
$C^{(K)}$ and that
\begin{eqnarray}
\mbox{dim}(Q^{(K)}) - \mbox{dim}(C^{(K)})&=&K(K-1)-(\mbox{dim}(U(K))- \mbox{dim}(USp(K)))\nonumber\\
&=&K(K-1)-K^2+\frac{1}{2}K(K+1)\nonumber\\
&=&\frac{1}{2}K(K-1),
\end{eqnarray}
we reach
\begin{equation} \label{J0}
J_{N,0} = e^{- \frac{NK}{2}} (2 \pi N)^{-\frac{K(K-1)}{4}}
\int_{C^{(K)}}\mu(dw) e^{Tr\left(Z^{\dagger}X Z X\right)} (1+O(N^{-1})).
\end{equation}
Collecting together (\ref{app:md}), (\ref{app:int3}) and (\ref{J0}) we find
\bea
\tilde{\rho}^{(N)}(x_1,x_2,\ldots, x_K )  = c_2(N,K)
\Delta(\mathbf{x}) \prod_{k=1}^{K} e^{-x_k^2} \int_{C^{(K)}}\mu(dw) e^{Tr\left(Z^{\dagger}X Z X\right)}
\left(1+ o(1)\right),\nonumber\\
\eea
where
\begin{eqnarray}
c_2(N,K) & = &  C(K)  \prod_{k=1}^{K} \left( |S_{N-k}|  \pi^{-\frac{N-k}{2}} \right)
\pi^{-\frac{K(K-1)}{2}} 2^{-\frac{(N-K)K}{2} } \\
&& \hspace{.2in} (N-K)^{\frac{(N-K)K}{2}}
(N-K)^{\frac{K(K-1)}{2}}
e^{- \frac{(N-K)K}{2}} (2 \pi (N-K))^{-\frac{K(K-1)}{4}}\nonumber
\end{eqnarray}
and $C(K)$ denotes a constant only depending on $K$. It is lengthy but straightforward
to check that $c_2(N,K) \to c_3(K) >0$ as $N \to \infty$
and hence that limiting modified density $\tilde{\rho}(x_1,x_2,\ldots, x_K ) = \lim_{N \to \infty}
\tilde{\rho}^{(N)}(x_1,x_2,\ldots, x_K ) $ exists and is given by
\begin{equation} \label{temp11}
\tilde{\rho}(x_1,x_2,\ldots, x_K )  = c_3(K)
\Delta(\mathbf{x}) \prod_{k=1}^{K} e^{-x_k^2} \int_{C^{(K)}}\mu(dw) e^{Tr\left(Z^{\dagger}X Z X\right)}.
\end{equation}
The explicit value of
$c_3(K)$ will be determined later by using properties of the densities $\tilde{\rho}$.
In the next subsection we will find a parameterisation
 of the integral in the right hand side of (\ref{temp11}), which will allow us to
 calculate it very efficiently using the standard tools  of random matrix theory.

\noindent
\textbf{Recasting as an integral over the unitary group.}
%%%%%%%%%%%%%%%%%%%%%%%%%%%%%%%%%%%%%%%%%%%%%%%%%%
An important property of the function $F$ is its invariance with respect to the following action of the unitary group
$U(K)$ on $Q^{(K)}$:
\begin{eqnarray}
U(K) \times Q^{(K)} & \to & U(K) \nonumber \\
(U,Z)  & \mapsto &  U Z U^T \in U(K). \label{app:actn}
\end{eqnarray}
Namely, for any $A \in U(K)$, we have $F(AZA^T)=F(Z)$.
The decomposition theorem for skew symmetric unitary matrices \cite{mehta1989matrix} states that
\bea\label{app:zdec}
Z=UJU^{T},
\eea
where $U$ is a unitary matrix, $J$ is the canonical symplectic matrix. Notice that (\ref{app:zdec}) does not
determine the unitary matrix $U$
uniquely: indeed $Z\rightarrow Z$ if $U\rightarrow US$, where $S$ is a unitary matrix satisfying $SJS^{T}=J$.
The set of such matrices is a subgroup of $U(K)$ called the symplectic group $USp(K)$, that is
\bea
USp(K)=\{S \in U(K)\mid SJS^T=J\}.
\eea
It can be checked that the critical manifold $C^{(K)}$ can be identified with the factor space of $U(K)$ with
respect to the action of $USp(K)$ on $U(K)$ via right multiplication, that is
\bea\label{app:cm}
C^{(K)} \cong U(K)/USp(K)
\eea
The $U(K)$-action (\ref{app:actn}) on $Q^{(K)}$ preserves the critical manifold and induces the $U(K)$-action on $C$.
Using the parameterisation (\ref{app:zdec}) of $C^{(K)}$ this induced action can be written explicitly:
\bea
U(K)\times C^{(K)} &\rightarrow& C^{(K)}, \nonumber \\
(A,[U]) &\mapsto& [AU], \label{app:actn1}
\eea
where $[U]$ is an equivalence class of $U \in U(K)$ with respect to right multiplications by elements of $USp(K)\subset U(K)$.
In the vicinity of a critical point $Z_c \in C^{(K)}$,
\bea\label{app:exp}
F(Z_c+\delta Z)=K+\frac{1}{2}Tr(\delta Z \dZ_c+Z_c\delta \dZ)^2+\ldots
\eea
We notice that the quadratic form describing the second order term in the above Taylor expansion of $F$ is $U(K)$-invariant and has
the maximal possible rank equal to $\frac{1}{2}K(K-1) =
\mbox{dim}(Q^{(K)})- \mbox{dim} (C^{(K)})$.
We rewrite the integral from (\ref{temp11}), using the mapping (\ref{app:cm}), as
\begin{equation}
\int_{C^{(K)}}\mu(dw) e^{Tr\left(Z^{\dagger}X Z X\right)}
=
\int_{U(K) /USp(K)} \hat{\mu}(dU) e^{Tr\left(Z(U)^{\dagger}X Z(U)X\right)} \label{app:io2})
\end{equation}
where $Z(U)$ is given by (\ref{app:zdec}) and $\hat{\mu}(dU)$ is the pull back of the measure $\mu$ on the critical manifold.
We can work out an explicit expression for $\mu$ in local co-ordinates on $C^{(K)}$ using the
general formula (\ref{app:meas}). We will not do this; instead we will characterise
$\hat{\mu}$ up to a multiplicative constant by establishing its symmetry with respect to the $U(N)$-action on $C^{(K)}$.
Recall that the measure $\mu$
is determined by the Lebesgue measure on $Q^{(K)}$ and the determinant
of the quadratic form in the right hand side of (\ref{app:exp}). It is easy to check
\begin{enumerate}
\item[(i.)] The Lebesgue measure $\Lambda$ and the quadratic form $Tr(\delta Z \dZ+Z\delta \dZ)^2$ on $Q^{(K)}$
are invariant with respect to the $U(K)$-action (\ref{app:actn}).
\item[(ii.)] The critical manifold $C^{(K)}$ is invariant with respect to the $U(K)$-action.
\item[(iii.)] The restriction of the quadratic form $Tr(\delta Z \dZ+Z\delta \dZ)^2$ on $Q^{(K)}$ to $C^{(K)}$ has maximal rank.
\end{enumerate}
A calculation employing elementary tools of differential geometry \cite{dubrovin1991modern} shows that
the above three observations imply the invariance of the measure
$\hat{\mu}$ with respect to the induced action  of $U(K)$  on the critical manifold $U(K)/USp(K)$ defined by (\ref{app:actn1}).
Therefore $\hat{\mu}$ is a  Haar measure on the symmetric space $U(K)/USp(K)$, which is unique up to normalization.
It is generally easier to work with integrals over the whole unitary group rather than a factor space. As we have established
already, the measure $\hat{\mu}$ is invariant with respect to the action of $U(K)$ on $C^{(K)}$.
Note that the function $Tr\left(Z(U)^{\dagger}X Z(U)X\right)$
which determines the integrand of (\ref{app:io2}) is also $U(K)$-invariant. Therefore, by Weyl's
theorem \cite{helgason1962differential} Chapter $X$,
\bea\label{app:io3}
\int_{U(K) /USp(K)} \hat{\mu}(dU) e^{Tr\left(Z(U)^{\dagger}X Z(U)X\right)}
=
\int_{U(K)}\mu_H(dU) e^{Tr\left(Z(U)^{\dagger}X Z(U)X\right)}
\eea
where $\mu_H$ is an appropriately normalized Haar measure on the unitary group. We will determine the
normalization factor later using the properties
of spin-spin correlation functions. Substituting (\ref{app:zdec}) into (\ref{app:io3}) we find that
\bea\label{app:io4}
\int_{U(K)}\mu_H(dU) e^{Tr\left(Z(U)^{\dagger}X Z(U)X\right)}
= \int_{U(K)}\mu_H(dU) e^{-Tr\left(JHJH^T\right)}
\eea
where $H$ is a Hermitian matrix with eigenvalues $x_1, x_2,\ldots, x_K$
given by $H=UXU^{\dagger}$.
Tracing back we find that the large $N$ behaviour of
$ \mathbb{E} \left[\prod_{m=1}^K\det\left(M_1-x_m\right)\right]$
 - the expected value of the product of $K$ characteristic polynomials in the real Ginibre ensemble -
turns out to be determined by the integration
of the symplectic-invariant Gaussian weight $\exp (-Tr JHJH^T)$ with respect to
unitary degrees of freedom. See \cite{sommers2007symplectic} for a discussion of the origin of the connection
between real Ginibre and symplectic ensembles. Let us note also that the integral (\ref{app:io4}) corresponds to $U(K)\rightarrow USp(K)$
symmetry breaking pattern according to the classification of Martin Zirnbauer \cite{mz3}.

We may rewrite the integral in (\ref{app:io4}) as
\begin{equation}
\label{app:f1}
\int_{U(K)}\mu_H(dU)
e^{Tr\left(HH^R\right)}
\end{equation}
where $H=UXU^{\dagger}$ and $H^R=JH^TJ^T$
is a 'symplectic' involution on the space of complex $K\times K$ matrices, see \cite{mehta2004random} for details.
Following Mehta, we will call matrix $M$ self-dual if $M=M^R$ and anti-self-dual if $M=-M^R$. It is easy to check
that any even-dimensional matrix can be uniquely represented as a sum of a self-dual and anti-self-dual matrices.
Let $ASD^{(K)}$ be the linear space of all anti-self dual $K\times K$ matrices. 
Combining (\ref{app:io4}), (\ref{app:f1}) with (\ref{temp11}) {\bf we obtain
the statement of Theorem \ref{thm3}}.

In order to derive the Pfaffian representation for the correlation functions,
we need to  perform the
 integration over the unitary group in the right hand side of (\ref{app:f1}). In order to achieve that, we use the following transformation found in
 \cite{mehtapandey1983},\cite{mehta1999some}:
\bea\label{app:mt}
e^{Tr H H^R}=Z_K e^{Tr H^2} \int_{ASD^{(K)}} \Lambda(dA) e^{Tr A^2+2\sqrt{2} Tr HA},
\eea
where $\Lambda(dA)$ is a Lebesgue measure on $ASD^{(K)}$, $Z_K$ is a normalization constant.
The absolute convergence of the above integral
can be checked using the decomposition theorem for anti-self-dual matrices,
see \cite{mehta1989matrix}: for any $A\in ASD^{(K)}$ there exists $V \in U(K)$ and  a diagonal
matrix $\Theta$ with entries $ \pm \theta_1, \pm \theta_2,\ldots
\pm \theta_{\frac{K}{2}}$, where $\theta_1\geq 0, \theta_2\geq 0, \ldots \theta_{K/2}\geq 0$, so that
\bea
A=iV\Theta V^{\dagger}.
\eea
Substituting (\ref{app:mt}) into (\ref{app:f1}) and using the invariance property of the Haar measure $\mu_H$ we get:
\bea
\int_{U(K)} \! \mu_H(dU) e^{Tr\left(HH^R\right)}
= C_K e^{\sum x_k^2} \int_{U(K)}\mu_H(dU) \int_{\mathbf{R}^{K/2}_{+}} \nu(d \theta)
e^{-Tr\left(\Theta^2\right)}e^{i2\sqrt{2}TrUXU^{\dagger}\Theta}\nonumber\\
\eea
where $\nu(d \theta)$ is the measure on the eigenvalues of unitary matrices induced by the marginalization
of the Lebesgue measure on $ASD^{(K)}$ over the unitary degrees of freedom, explicitly \cite{mehta2004random}
\bea
\nu(d \theta)=\Delta\left(\pm \theta_1, \pm \theta_2,\ldots
\pm \theta_{\frac{K}{2}}\right) \prod_{k=1}^{K/2} \theta_k \, d \theta_k.
\eea
We now allow the constants $C_K$, depending only on $K$, to change value from line to line.
We can integrate over the unitary group $U(K)$ using Harish-Chandra-Itzykson-Zuber
formula \cite{harish1957differential}, \cite{itzykson1980planar}. The result is
\begin{eqnarray}
&& \hspace{-.3in} \int_{U(K)}\mu_H(dU) e^{Tr\left(HH^R\right)} \\
& = &
C_K e^{\sum x_k^2} \Delta(x)^{-1}
 \int_{\mathbf{R}^{K/2}_{+}} \prod_{k=1}^{K/2} \theta_k e^{-\theta_k^2} \, d\theta_k \;
det \left[e^{i2\sqrt{2}x_j\Theta_{ll}}: 1\leq j,l \leq K \right].\nonumber
\end{eqnarray}
The remaining integration over the singular values $\theta_{1}, \theta_{2}, \ldots, \theta_{K/2}$ is carried out
using the de Bruijn formula \cite{de1955some}, yielding
\bea
\int_{U(K)}\mu_H(dU) e^{Tr\left(HH^R\right)}
= C_K e^{\sum x_k^2} \Delta(x)^{-1} Pf\left[(x_i-x_j)e^{-2(x_i-x_j)^2}\right]_{1\leq i,j\leq K}.
\eea

Combined with (\ref{temp11}) this yields
\begin{equation} \label{temp12}
\tilde{\rho}(x_1,x_2,\ldots, x_K )  = C_k  Pf\left[(x_i-x_j)e^{-2(x_i-x_j)^2}\right]_{1\leq i,j\leq K}.
\end{equation}
The correlation functions of spins can be formally computed by integrating $\tilde{\rho}$ with respect to space variables:
\begin{eqnarray}\label{eq:ss1d}
\mathbb{E}\left(\prod_{k=1}^K s_{x_k}(M_1) \right) = (-2)^K\left(\prod_{k=1}^K \int_{-\infty}^{x_k} dy_k\right)
\tilde{\rho}(y_1,y_2,\ldots, y_{k}).
\end{eqnarray}
This leads to the spin-spin correlation function
\begin{equation} \label{temp13}
\mathbb{E} \left(\prod_{k=1}^K s_{k_k}(M_1) \right)=
C_k  Pf\left[ \int_{x_i-x_j} e^{-2 z^2} dz \right]_{1\leq i,j\leq K}.
\end{equation}
The constants $C_k$ can be found inductively in $k$ by allowing $x_{2k} \downarrow x_{2k-1}$, and
noting that $\rho(x_1,x_1) =1$. This yields $C_K = (8/\pi)^{K/4}$.
Expression (\ref{temp13}) coincides with the single time spin-spin correlation
for a system of one dimensional annihilating Brownian motions under the maximal entrance law, see \cite{EJP942}
for details. As shown in\cite{EJP942},
substituting (\ref{temp13}) into formula (\ref{eq:dens})
we get the first statement of Corollary $9$ of \cite{borodin2009ginibre}.
%%%%%%%%%%%%%%%%%%%%%%

\noindent
{\bf Acknowledgement.} We are grateful to Freddy Bouchet, Neil O'Connell, Krzysztof Gawedzki,
Alice Guionnet, John Rawnsley and Stefano Ruffo for helpful discussions.

%%%%%%%%%%%%%%%%%%%%%%%%%%%%%%%%%%%%%%%%%%%%%%%%%%%%%%%%%%%%%%%%%%%%5
%\bibliography{rg12nov_final}{}
%\bibliographystyle{plain}
%%%%%%%%%%%%%%%%%%%%%%%%%%%%%%%%%%%%%%%%%%%%%%%%%%%%%%%%%%%%%%%%%%%%%

\end{document}